# DOES MEDIAN FILTERING TRULY PRESERVE EDGES BETTER THAN LINEAR FILTERING?

By Ery Arias-Castro and David L. Donoho

*University of California, San Diego and Stanford University*

Image processing researchers commonly assert that "median filtering is better than linear filtering for removing noise in the presence of edges." Using a straightforward large-$n$ decision-theory framework, this folk-theorem is seen to be false in general. We show that median filtering and linear filtering have similar asymptotic worst-case mean-squared error (MSE) when the signal-to-noise ratio (SNR) is of order 1, which corresponds to the case of constant per-pixel noise level in a digital signal. To see dramatic benefits of median smoothing in an asymptotic setting, the per-pixel noise level should tend to zero (i.e., SNR should grow very large).

We show that a two-stage median filtering using two very different window widths can dramatically outperform traditional linear and median filtering in settings where the underlying object has edges. In this two-stage procedure, the first pass, at a fine scale, aims at increasing the SNR. The second pass, at a coarser scale, correctly exploits the nonlinearity of the median.

Image processing methods based on nonlinear partial differential equations (PDEs) are often said to improve on linear filtering in the presence of edges. Such methods seem difficult to analyze rigorously in a decision-theoretic framework. A popular example is mean curvature motion (MCM), which is formally a kind of iterated median filtering. Our results on iterated median filtering suggest that some PDE-based methods are candidates to rigorously outperform linear filtering in an asymptotic framework.

## 1. Introduction.

1.1. *Two folk theorems.* Linear filtering is fundamental for signal processing, where often it is used to suppress noise while preserving slowly varying signal. In image processing, noise suppression is also an important task; however, there have been continuing objections to linear filtering of images since at least the 1970s, owing to the fact that images have edges.

In its simplest form, linear filtering consists of taking the average over a sliding window of fixed size. Indeed, linear filtering with fixed window size $h$ "blurs out" the edges, causing a bias of order $O(1)$ in a region of width $h$ around edges. This blurring can be visually annoying and can dominate the mean-squared error.

Median filtering—taking the median over a sliding window of fixed size— was discussed already in the 1970s as a potential improvement on linear filtering in the "edgy" case, with early work of Matheron and Serra on morphological filters [24, 33] in image analysis and (in the case of 1-d signals) by Tukey and collaborators [22, 37, 38].

To this day simple median filtering is commonly said to improve on linear filtering in "edgy" settings [2, 8, 14, 16, 36]—such a claim currently appears in the *Wikipedia* article on median filtering [1]. Formally, we have the

MEDIAN FOLK THEOREM.   *Median filtering outperforms linear filtering for suppressing noise in images with edges.*

Since the late 1980s, concern for the drawbacks of linear filtering of images with edges has led to increasingly sophisticated proposals. In particular, inspired by seminal work of Mumford and Shah [26] and Perona and Malik [27], a whole community in applied mathematics has arisen around the use of nonlinear partial-differential equations (PDEs) for image processing— including noise suppression [25, 31, 34].

A commonly heard claim at conferences in image processing and in applied mathematics boils down to the following:

PDE FOLK THEOREM.   *PDE-based methods outperform linear filtering for suppressing noise in images with edges.*

1.2. *A challenge to asymptotic decision-theory.*   While these folk theorems have many believers, they implicitly pose a challenge to mathematical statisticians.

Linear filtering of the type used in signal and image processing has also been of interest to mathematical statisticians in implicit form for several decades. Indeed, much of nonparametric regression, probability density estimation and spectral density estimation is in some sense carried out with kernel methods—a kind of "linear filter"—and there is extensive literature documenting the optimality of such linear procedures in certain cases. In many cases, the correspondence of the underlying bias-variance analysis with the kind of analyses being done in signal processing is quite evident.

During the last two decades, mathematical statisticians have succeeded in showing that for models of images with edges, nonlinear methods can indeed outperform linear ones in a minimax sense [5, 10, 11]. However, the nonlinear methods which have been analyzed fully rigorously in a decision-theoretic framework are somewhat different than the median and PDE cases above; examples include methods of wavelet shrinkage and other harmonic analysis techniques [5, 9, 10, 11].

So within the decision-theoretic framework it is rigorously possible to do better than linear filtering, but by methods somewhat different than those covered by the folk theorems above.

The great popularity of median- and PDE-based nonlinear filtering prompted us to evaluate their performance by a rigorous approach within the decision-theoretic framework of mathematical statistics. Three conclusions emerge:

- The Median-filtering folk theorem is *false* in general.
- In an apparently meaningless special case—where the noise level per pixel is negligible—the Median-filtering folk theorem is true.
- A modified notion of median filtering—applying two passes in a multiscale fashion—*does* improve on linear filtering, as we show here.

Before explaining these conclusions in more detail, we make a few remarks.

- Tukey's emphasis with median filters was always on iterating median filters sequentially applying medians over windows of different widths at different stages, as we do here. We believe that Tukey's intuition about the benefits of median filtering actually applied to this iterated form; but this intuition was either never formalized in print or has been forgotten with time.
- The iterated median scheme we are able to analyze in this paper simply involves two passes of medians at two very different scales.

Finally, we believe our results are part of a bigger picture:

- There is a formal connection between certain nonlinear PDEs used for image processing (i.e., Mean Curvature Motion and related PDEs) and iterated medians.
- Nonlinear PDE-based methods in the form usually proposed by applied mathematicians seem quite difficult to analyze within the decision-theoretic framework; this seems a looming challenge for mathematical statisticians.
- Our iterated median scheme seems related to such PDE-based methods. It is perhaps less elegant than full nonlinear PDEs but is rigorously analyzed here.
- Because of results reported here, we now suspect some subset of the PDE folk theorem may well be true.

1.3. *The framework.* We describe the framework below in any dimension $d$; however, this paper focuses on dimensions $d = 1$ (signal processing) and $d = 2$ (image processing).

Consider the classical problem of recovering a function $f : [0,1]^d \mapsto [0,1]$ from equispaced samples corrupted by additive white noise. Here we observe $Y_n = (Y_n(\mathbf{i}))$ defined by

$$(1.1) \qquad Y_n(\mathbf{i}) = f(\mathbf{i}/n) + \sigma Z_n(\mathbf{i}), \qquad \mathbf{i} \in \mathbb{I}_n^d,$$

where $\mathbb{I}_n = \{1, \ldots, n\}$, $d$ is the dimension of the problem (here, $d = 1$ or 2), $\sigma > 0$ is the noise level and $Z_n$ is white noise with distribution $\Psi$. Only $Y_n$, $\sigma$ and $\Psi$ are known. Though unknown, $f$ is restricted to belong to some class $\mathcal{F}$ of functions over $[0,1]^d$ with values in $[0,1]$; several such $\mathcal{F}$ will be explicitly defined in Sections 2, 3 and 4.

By linear filtering we mean the following variant of moving average. Fix a window size $h \geq 0$ and put

$$L_h[Y_n](\mathbf{i}) = \operatorname{Average}\{Y_n(\mathbf{j}) : \mathbf{j} \in \mathcal{W}[n, h](\mathbf{i})\};$$

here $\mathcal{W}[n, h](\mathbf{i})$ denotes the discrete window of radius $nh$ centered at $\mathbf{i} \in \mathbb{I}_n^d$:

$$\mathcal{W}[n, h](\mathbf{i}) = \{\mathbf{j} \in \mathbb{I}_n^d : \|\mathbf{j} - \mathbf{i}\| \leq nh\}.$$

Similarly, by median filtering we mean

$$M_h[Y_n](\mathbf{i}) = \operatorname{Median}\{Y_n(\mathbf{j}) : \mathbf{j} \in \mathcal{W}[n, h](\mathbf{i})\}.$$

(More general linear and median filters with general kernels add no dramatically new phenomena in settings of interest to us, i.e., where signals are discontinuous, so we ignore them.)

Following a traditional approach in mathematical statistics [21], the performance of an estimator $T[Y_n]$ is measured according to its worst-case risk over the functional class of interest $\mathcal{F}$ with respect to mean-squared error (MSE):

$$\mathcal{R}_n(T; \mathcal{F}) = \sup_{f \in \mathcal{F}} \mathcal{R}_n(T; f),$$

where

$$\mathcal{R}_n(T; f) = \frac{1}{n^d} \sum_{\mathbf{i} \in \mathbb{I}_n^d} \mathbb{E}[(T[Y_n](\mathbf{i}) - f(\mathbf{i}/n))^2].$$

We consider for $\mathcal{F}$ certain classes of piecewise Lipschitz functions. When the noise distribution $\Psi$ is sufficiently nice—Gaussian, for example—we show that linear filtering and median filtering have worst-case risks with the same rates of convergence to zero as $n \to \infty$. This contradicts the Median folk theorem.

Our conclusion does not rely on misbehavior at any farfetched function $f \in \mathcal{F}$. In dimension $d = 1$, linear filtering and median filtering exhibit the same worst-case rate of convergence already at the simple step function $f(x) = 1_{\{x > 1/2\}}$—the simplest model of an edge.

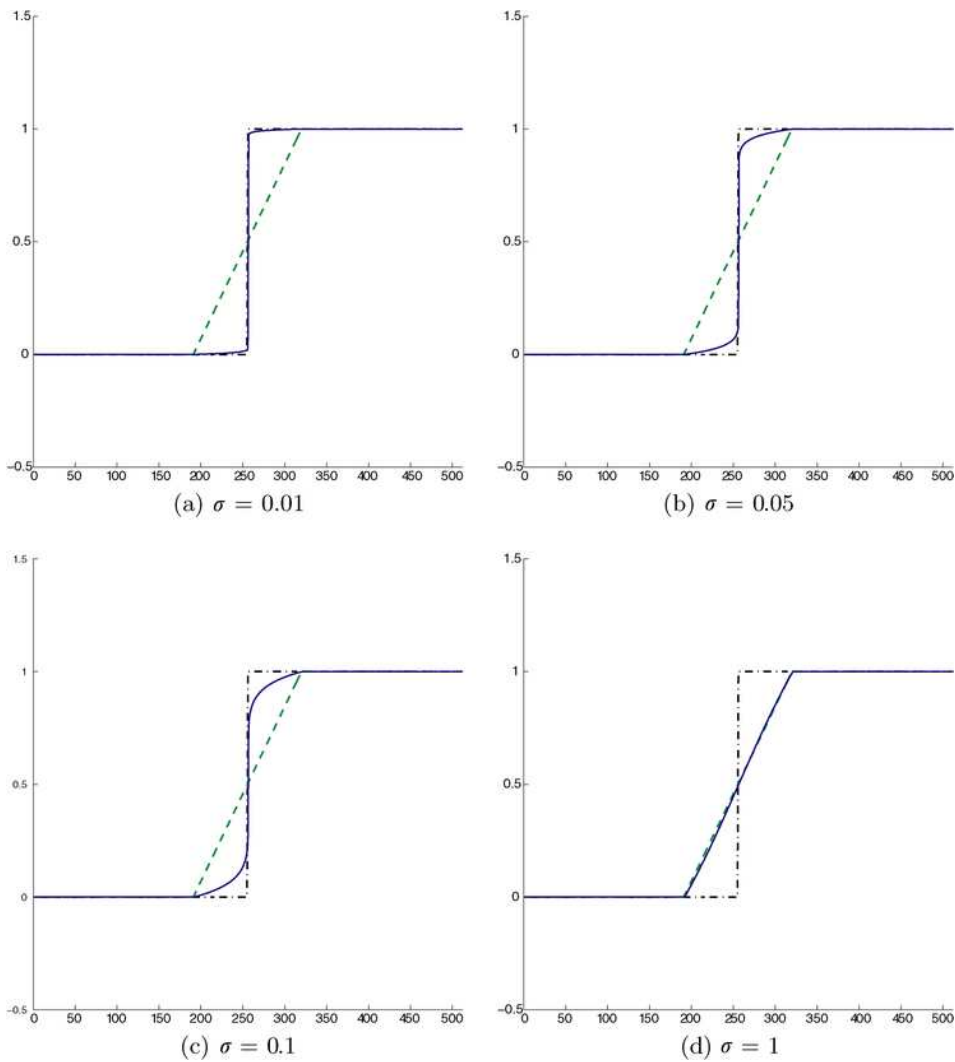

Fig. 1. *The dash-dotted (black) line represents the noiseless object $f = 1_{[1/2,1]}$; the dashed (green) line represents the expected value of linear filtering; the solid (blue) line represents the expected value of median filtering. The noise is Gaussian, the smoothing window size is $h = 0.125$ and the sample size $n = 512$. Only at very small $\sigma$ is the bias of the median qualitatively superior to the bias of linear filtering.*

1.4. *The underlying phenomenon.* The misbehavior of median filtering can be traced to the fact that, *for a signal-to-noise ratio of order 1 (specifically, for $\sigma = 1$),* its bias is of order 1 in a region of width $h$ near edges; this behavior is virtually identical to that of linear filtering. Figure 1 illustrates this situation in panel (d).

However, the figure also illustrates, in panels (a), (b), (c), another phenomenon: for very low noise levels $\sigma$, that is, very high signal-to-noise ratios, the bias of the median behaves dramatically differently than the bias of linear filtering. In particular, the bias is not large over an interval comparable to the window width, but only over a much smaller interval. In fact, as $\sigma \to 0$, the bias vanishes away from the edge. We call this the:

TRUE HOPE OF MEDIAN FILTERING. At very low noise levels, median filtering can dramatically outperform linear filtering.

At first glance this seems utterly useless: why should we care to remove noise when there is almost no noise? On reflection, a useful idea emerges. Suppose we filter in stages, at the first stage using a relatively narrow window width—much narrower than we would ordinarily use in a one-stage process—and at the second stage using a somewhat wider window width. The result may well achieve the noise reduction of the combined two-stage smoothing with much smaller bias near edges.

HEURISTIC OF ITERATED MEDIAN FILTERING. Iterated median filtering, in which the data are first median-filtered lightly, at a fine scale, followed by a coarse-scale median filter, may outperform linear filtering.

Note that the same idea, applied to linear filtering, would achieve little. The composition of two linear filters can always be achieved by a single linear filter with appropriate kernel. And such weights do not change the qualitative effect of edges.

TABLE 1

**Summary of results.** *Rates of convergence to zero of worst-case MSE, and of optimal window width, for different methods in dimensions $d = 1,\ 2$*

| Dimension | $d = 1$ | | $d = 2$ | |
|---|---|---|---|---|
| Technique | Rate | Window width | Rate | Window width |
| Linear filter | $n^{-1/2}$ | $n^{-1/2}$ | $n^{-2/3}$ | $n^{-2/3}$ |
| Median filter | $n^{-1/2}$ | $n^{-1/2}$ | $n^{-2/3}$ | $n^{-2/3}$ |
| Two-scale median filter | $n^{-2/3}$ | $n^{-2/3},\ n^{-1/3}$ | $n^{-6/7}$ | $n^{-6/7},\ n^{-4/7}$ |
| Edge-free optimal | $n^{-2/3}$ | | $n^{-1}$ | |

In each case, the underlying class of functions has Lipschitz smoothness away from edges. Results compiled from theorems below. Note: here $n$ is the signal width in pixels, not the sample size. The sample size is $n^d$ in dimension $d$.

1.5. *Results of this paper.* Table 1 compiles the worst-case risk rates for linear filtering, median filtering and two-scale median filtering over piecewise Lipschitz function classes. The last line displays the minimax rates for Lipschitz functions *without discontinuities*. Note that the rates would be the same for classes of functions with higher degree of smoothness away from the edges—indeed, the three filtering methods considered here have worst-case risk of same order of magnitude as their MSE for a simple step function such as $f(x) = 1_{\{x > 1/2\}}$ in dimension $d = 1$.

Notice also that, in dimension $d = 1$ our two-scale median filtering achieves the minimax rate for edge-free Lipschitz function classes. This is not the case in dimension $d = 2$, where other methods are superior [5, 21].

1.6. *Iterated medians.* As mentioned above, in the late 1960s Tukey already proposed the use of iterated medians, although the motivations remained unclear to many at the time. In his proposals, different scales were involved at different iterations, although the scales were relatively similar from the current viewpoint.

Significantly, iterated median filtering converges in some sense to Mean Curvature Motion (MCM), a popular PDE-based technique. In fact, there is a wider link connecting PDE-based methods and iterative filtering; in particular, iterated linear filtering converges to the Heat equation. See, for example, [7, 8, 15]. Although MCM is highly nonlinear and hard to analyze, the heuristic above gives a hint that MCM might improve on linear filtering.

1.7. *Prior literature.* Several papers analyze the performance of median filtering numerically using simulations. For example, in [17], the authors derive exact formulas for the distribution of the result of applying median filtering to a simple noisy edge like $f_0$, and use computer-intensive simulations to provide numerical values. A similar approach is found in [19].

Closer to the present paper, [20] compares linear filtering and median filtering in the context of *smooth* functions and shows that they have minimax rates of same order of magnitude. We will show here that the same holds for functions with discontinuities.

An extensive, but unrelated, body of literature explores the median's ability to suppress outliers, for example, [13] and also [28], which consider the case of smoothing a one-dimensional signal corrupted with impulsive noise. Using a similar framework, Donoho and Yu [12] study a pyramidal median transform.

1.8. *Contents.* In Sections 2 and 3, we consider one-dimensional and two-dimensional signals, respectively, in the constant-noise level case. In Section 4 we consider per-pixel noise level tending to 0. In Section 5, we introduce a two-scale median filtering and formulate results quantifying its performance. The proofs are postponed to the latter sections.

Though our results can be generalized to higher-dimensional signals and other smoothness classes, and can accommodate more sophisticated kernels, we choose not to pursue such extensions and generalizations here.

1.9. *Notation.* Below, we denote comparable asymptotic behavior of sequences $(a_n), (b_n) \in \mathbb{R}^{\mathbb{N}}$ using $a_n \asymp b_n$, meaning that the ratio $a_n/b_n$ is bounded away from 0 and $\infty$ as $n$ becomes large.

**2. Linear filtering and median filtering in dimension 1.** Consider the model (1.1) introduced in the Introduction for the case of dimension 1 ($d = 1$). We now explicitly define the smoothness class of interest.

DEFINITION 2.1. The *local Lipschitz constant* of the function $f$ at $x$ is

$$\mathrm{Lip}_x(f) = \limsup_{\varepsilon \to 0} \sup_{|y-x| \le \varepsilon, y \neq x} \frac{|f(y) - f(x)|}{|y - x|}.$$

A function $f : [0,1] \mapsto \mathbf{R}$ will be called *essentially local Lipschitz* if the essential supremum of the local Lipschitz constant on $[0,1]$ is finite.

The function $x 1_{\{x > 1/2\}}(x)$ is essentially local Lipschitz but not Lipschitz; it has a local Lipschitz constant $\le 1$ almost everywhere on $[0,1]$, but jumps across the line $x = 1/2$. More generally, piecewise polynomials without continuity constraints at piece boundaries are essentially local Lipschitz and yet neither Lipschitz nor continuous.

DEFINITION 2.2. Fix $N \ge 1$ and $\beta > 0$. The class of *punctuated-Lipschitz* functions $\mathrm{PLIP} = \mathrm{PLIP}(\beta, N)$ is the collection of functions $f : [0,1] \mapsto [0,1]$ with local Lipschitz constant bounded by $\beta$ on the complement of some finite set $(x_i)_{i=1}^{N} \subset (0,1)$.

THEOREM 2.1. *Assume $\Psi$ has mean 0 and finite variance. Then,*

$$\inf_{h > 0} \mathcal{R}_n(L_h; \mathrm{PLIP}) \asymp n^{-1/2}, \qquad n \to \infty.$$

This result can be proven using standard bias-variance trade-off ideas, and needs only simple technical ingredients such as uniform bounds on functions and derivatives. We explain in Section 7 that it can be inferred from existing results, but then proceed to give a proof; this proof sets up a bias-variance trade-off framework suitable for several less elementary situations which come later.

To analyze the behavior of median filtering, we must obtain uniform bounds on the stochastic behavior of empirical quantiles; these are laid out in Section 15 below. To enable such bounds, we make the following assumptions on the noise distribution $\Psi$:

[Shape] $\Psi$ *has density* $\psi$ *with respect to the Lebesgue measure on* $\mathbb{R}$, *with* $\psi$ *unimodal, continuous and symmetric about* $0$.
[Decay] $\zeta(\Psi) := \sup\{s > 0 : \psi(x)(|x| + 1)^s$ *is bounded*$\} > 1$.

Note that the normal, double-exponential, Cauchy and uniform distributions centered at 0 all satisfy both [Shape] and [Decay]. These conditions permit an efficient proof of the following result—see Section 8; the conditions could probably be relaxed considerably, leading to a more difficult proof.

THEOREM 2.2. *Assume* $\Psi$ *satisfies* [Shape] *and* [Decay]. *Then,*
$$\inf_{h > 0} \mathcal{R}_n(M_h; \text{PLIP}) \asymp n^{-1/2}, \qquad n \to \infty.$$

For piecewise Lipschitz functions corrupted with white additive Gaussian noise (say) having constant per-pixel noise level, Theorems 2.1 and 2.2 show that linear filtering and median filtering have risks of same order of magnitude $\mathcal{R}_n \asymp n^{-1/2}$—the same holds for any noise distribution $\Psi$ with finite variance satisfying [Shape] and [Decay]. In that sense, the Median folk theorem is contradicted.

The proof shows that the optimal order of magnitude for the width $h_n$ of the median smoothing window obeys $h_n \asymp n^{-1/2}$—again the same as in linear filtering.

## 3. Linear filtering and median filtering in dimension 2.

Consider model (1.1) of the Introduction in the case of dimension 2 ($d = 2$); our "signals" are now digital images. Our smoothness class here is a class of *cartoon images*, which are piecewise functions that are smooth except for discontinuities along smooth curves—see also [6, 9, 21].

Just as in $d = 1$, we have the notion of local Lipschitz constant. In $d = 2$, the function $x1_{\{x > 1/2\}}(x, y)$ is essentially local Lipschitz but not Lipschitz; this has a local Lipschitz constant bounded by 1 almost everywhere, but the function jumps as we cross the line $x = 1/2$. More generally, cartoon images have the same character: essentially local Lipschitz and yet neither Lipschitz nor continuous. Such cartoon images, of course, have jumps along collections of regular curves; we formalize such collections as follows.

DEFINITION 3.1. A finite collection of rectifiable planar curves will be called a *complex*. Fix $\lambda > 0$ and let $\Gamma = \Gamma(\lambda)$ denote the class of rectifiable curves in $[0, 1]^2$ with length at most $\lambda$. Let $\mathcal{C}(N, \lambda)$ denote the collection of complexes composed of at most $N$ curves from $\Gamma(\lambda)$.

DEFINITION 3.2. Fix $N \geq 1$ and $\beta > 0$. The class of *curve-punctuated Lipschitz* functions $\text{CPLIP} = \text{CPLIP}(\lambda, \beta, N)$ is the collection of functions $f : [0, 1]^2 \mapsto [0, 1]$ having local Lipschitz constant bounded by $\beta$ on the complement of a $\mathcal{C}(N, \lambda)$-complex.

Informally, such functions are "locally Lipschitz away from edges" and indeed can be viewed as models of "cartoons." We prove the next two results in Sections 10 and 11, respectively.

THEOREM 3.1. *Assume $\Psi$ has mean 0 and variance 1. Then*

$$\inf_{h>0} \mathcal{R}_n(L_h; \text{CPLIP}) \asymp n^{-2/3}, \qquad n \to \infty.$$

THEOREM 3.2. *Assume $\Psi$ satisfies* [Shape] *and* [Decay]. *Then,*

$$\inf_{h>0} \mathcal{R}_n(M_h; \text{CPLIP}) \asymp n^{-2/3}, \qquad n \to \infty.$$

The situation parallels the one-dimensional case. In words, for cartoon images corrupted with white additive noise with constant per-pixel noise level, linear filtering and median filtering have risks of same order of magnitude $\mathcal{R}_n \asymp n^{-2/3}$. Again, the Median folk theorem is contradicted.

The proofs show that the width of the optimal smoothing window for either type of smoothing is $\asymp n^{-2/3}$.

**4. Linear filtering and median filtering with negligible per-pixel noise level.** The analysis so far assumes that the noise level is comparable to the signal level.

For very low-noise-per-pixel level and discontinuities well-separated from the boundary and each other, the situation is completely different: *the Median folk theorem holds true*.

*Preliminary remark.* We will see in Sections 9 and 12 that, both in dimensions $d = 1$ and $d = 2$, linear filtering does not improve on no-smoothing if $\sigma_n n^{1/2} = O(1)$, while median filtering improves on no-smoothing if $\sigma_n n \to \infty$. In Theorems 4.1 and 4.2 below we therefore exclude the situation $\sigma_n n = O(1)$.

DEFINITION 4.1. The finite point set $(x_i)_{i=1}^N \subset [0, 1]$ is called *well-separated* with separation constant $\eta > 0$ if (i) each point is at least $\eta$-separated from the boundary $\{0, 1\}$:

$$\min(x_i, 1 - x_i) \geq \eta \qquad \forall i$$

and (ii) each point is at least $\eta$-separated from every other point:

$$|x_j - x_i| \geq \eta \qquad \forall i, j.$$

DEFINITION 4.2. Let SEP-PLIP = SEP-PLIP$(\eta, \beta, N)$ denote the class of functions in PLIP$(\beta, N)$ which have local Lipschitz constant $\leq \beta$ on the complement of some $\eta$-well-separated set $(x_i)_{i=1}^N$.

THEOREM 4.1. *Let $\Psi$ satisfy* [Shape] *and* [Decay] *and let the per-pixel noise level tend to zero with increasing sample size, $\sigma = \sigma_n \to 0$ as $n \to \infty$, with $\sigma_n n \to \infty$. Then,*

$$\inf_{h>0} \mathcal{R}_n(M_h; \text{SEP-PLIP}) = o\left(\inf_{h>0} \mathcal{R}_n(L_h; \text{SEP-PLIP})\right), \qquad n \to \infty.$$

The proof is in Section 9, where we provide more explicit bounds for the risks of linear and median filtering.

In dimension $d = 2$, we again have that for negligible per-pixel noise level the Median folk theorem holds true. To show this, we need the hypothesis that the discontinuity curves are well-separated from the boundary, and from each other.

DEFINITION 4.3 (Well-separated complex). Let $d(A, B)$ denote Hausdorff distance between compact sets $A$ and $B$. A complex $C = (\gamma_i)$ of rectifiable curves in $[0,1]^2$ is said to be *well-separated* with separation parameter $\eta > 0$, if (i) the curves are separated from the boundary of the square:

$$d(\gamma_i, \text{BDRY}[0,1]^2) \geq \eta \qquad \forall i$$

and (ii) the curves are separated from each other:

$$d(\gamma_i, \gamma_j) \geq \eta \qquad \forall i, j.$$

We also need that the curves are well-separated from themselves (i.e., do not loop back on themselves). Formally, we need the condition

DEFINITION 4.4 ($C^2$ Chord-arc curves). Fix parameters $\lambda, \kappa, \theta$. Let $\Gamma_2 = \Gamma_2(\lambda, \kappa, \theta)$ be the collection of planar $C^2$ curves $\gamma$ with curvature bounded by $\kappa$ and chord-arc ratio bounded by $\theta$:

$$\forall s < t \qquad \frac{t-s}{|\gamma(t) - \gamma(s)|} \leq \theta \quad \text{and} \quad \frac{\text{length}(\gamma) - t + s}{|\gamma(t) - \gamma(s)|} \leq \theta.$$

Related classes of curves appear in, for example, Section 5.3 of [21]. Note that curves with bounded chord-arc ratio appear in harmonic analysis related to potential theory, for example, [32].

DEFINITION 4.5. Let SEP-CPLIP = SEP-CPLIP$(\lambda, \theta, \kappa, \eta, \beta, N)$ be the collection of curve-punctuated Lipschitz functions with local Lipschitz constant bounded by $\beta$ on the complement of an $\eta$-well-separated $\mathcal{C}(N, \lambda)$-complex of $\kappa, \theta$ chord-arc curves.

For an example, let $f = 1_D$, where $D$ is the disk of radius $1/4$ centered at $(1/2, 1/2)$. Then the exceptional complex $\mathcal{C} = \{\gamma_1\}$, where

$$\gamma_1(\theta) = (\tfrac{1}{2}, \tfrac{1}{2})' + \tfrac{1}{4}(\cos(\theta), \sin(\theta))', \qquad \theta \in [0, 2\pi).$$

For this example, $f \in$ SEP-CPLIP with parameters

$$N \geq 1, \qquad \lambda \geq \frac{\pi}{2}, \qquad \beta \geq 0, \qquad \theta \geq \frac{\pi}{2}, \qquad \kappa \geq 4, \qquad \eta \leq \tfrac{1}{4}.$$

As this example shows, we may choose parameters so that the classes $\Gamma$ and SEP-CPLIP are nonempty. In the sequel, we assume this has been done, so that SEP-CPLIP contains the $f$ just given.

THEOREM 4.2. *Assume $\Psi$ satisfies* [Shape] *and* [Decay] *and that the per-pixel noise level tends to zero with increasing sample size, $\sigma = \sigma_n \to 0$ as $n \to \infty$, with $\sigma_n n \to \infty$. Then,*

$$\inf_{h>0} \mathcal{R}_n(M_h; \text{SEP-CPLIP}) = o\left(\inf_{h>0} \mathcal{R}_n(L_h; \text{SEP-CPLIP})\right), \qquad n \to \infty.$$

The proof is in Section 12.

**5. Iterated and two-scale median filtering.** Iterated (or repeated) median filtering applies a series of median filters $M_{h_1,\dots,h_m}[Y_n] \equiv M_{h_m} \circ \cdots \circ M_{h_1}[Y_n]$. That is, median filtering with window size $h_1$ is first applied, then median filtering with window size $h_2$ is applied to the resulting signal, and so on. In the 1970s Tukey advocated such compositions of medians in connection with $d = 1$ signals, for example, applying medians of lengths 3, 5, and 7 in sequence—possibly, along with other operations, including linear filtering. Here we are interested in much longer windows than Tukey, in fact in windows that grow large as $n$ increases.

Tukey also advocated the iteration of medians until convergence—his so-called "3R" median filter applies running medians of three repeatedly until no change occurs. The mathematical study of repeated medians $M_{h_1,\dots,h_m}$ is a challenging endeavor, however, because of the strong dependency that median filtering introduces with every pass, though [3, 4] attempt to carry out just such studies, in the situation where there is only noise and no signal. See also Rousseeuw and Bassett [29].

Here, inspired by the intuition supplied in Figure 1 and by the results of the previous section, we consider two-scale median filtering. The first pass aims at increasing the signal-to-noise ratio, so the second pass can exploit the promising characteristics of median filtering at high signal-to-noise ratio.

We describe the process in dimension $d$. For $h > 0$, consider the squares

$$B_{\mathbf{k}}^h = [k_1 nh + 1, (k_1 + 1)nh] \times \cdots \times [k_d nh + 1, (k_d + 1)nh],$$

where $\mathbf{k} = (k_1, \dots, k_d) \in \mathbb{N}^d$ with $0 \leq k_j < 1/h$. Fix $0 < h_1 < h_2 < 1$ and define $M^{h_1, h_2}[Y_n]$ as follows:

1. For $\mathbf{k} = (k_1, \ldots, k_d) \in \mathbb{N}^d$ with $0 \le k_j < 1/h_1$, define

$$Y_n^{h_1}(\mathbf{k}) = \text{Median}\{Y_n(\mathbf{i}) : \mathbf{i} \in B_{\mathbf{k}}^{h_1}\};$$

thus $Y_n^{h_1}$ is a coarsened version of $Y_n$.

2. For $\mathbf{i} \in B_{\mathbf{k}}^{h_1}$, define

$$M^{h_1, h_2}[Y_n](\mathbf{i}) = M_{h_2}[Y_n^{h_1}](\mathbf{k}).$$

In words, we apply median filtering to the coarse-scale version $Y_n^{h_1}$ and get back the fine-scale version by crude piecewise interpolation. The following results, for $d = 1, 2$, respectively, are proven in Sections 13 and 14.

THEOREM 5.1.   *Assume* $\Psi$ *satisfies* [Shape] *and* [Decay]. *Then,*

$$\inf_{0 < h_1 < h_2} \mathcal{R}_n(M^{h_1, h_2}; \text{SEP-PLIP}) \asymp n^{-2/3}, \qquad n \to \infty.$$

THEOREM 5.2.   *Assume* $\Psi$ *satisfies* [Shape] *and* [Decay]. *Then,*

$$\inf_{0 < h_1 < h_2} \mathcal{R}_n(M^{h_1, h_2}; \text{SEP-CPLIP}) \asymp n^{-6/7}, \qquad n \to \infty.$$

Compare Theorems 5.1 and 5.2 with Theorems 2.2 and 3.2, respectively. In dimension $d = 1$, the rate improves from $O(n^{-1/2})$ to $O(n^{-2/3})$; in dimension $d = 2$, from $O(n^{-2/3})$ to $O(n^{-6/7})$. Hence, with carefully chosen window sizes, our two-scale median filtering outperforms linear filtering. The optimal choices are $h_1 \asymp n^{-2/3}$ and $h_2 \asymp n^{-1/3}$ for $d = 1$; $h_1 \asymp n^{-6/7}$ and $h_2 \asymp n^{-4/7}$ for $d = 2$.

Though this falls short of proving that indefinitely iterated median filtering of the sort envisioned by Tukey dramatically improves on linear filtering or that the PDE folk theorem is true for Mean Curvature Motion, it certainly suggests hypotheses for future research in those directions.

## 6. Tools for analysis of medians.
Before proceeding step-by-step with proofs of the theorems announced above, we isolate some special facts about medians which are used frequently and which ultimately drive our analysis.

6.1. *Elementary properties.*   Let $\text{Med}_n(\cdot)$ denote the empirical median of $n$ numbers. We make the following obvious but essential observations:

- *Monotonicity.* If $x_i \le y_i$, $i = 1, \ldots, n$,

(6.1)                    $$\text{Med}_n(x_1, \ldots, x_n) \le \text{Med}_n(y_1, \ldots, y_n).$$

- *Lipschitz mapping.*

(6.2)        $$|\text{Med}_n(x_1, \ldots, x_n) - \text{Med}_n(y_1, \ldots, y_n)| \le \max_{i=1}^{n} |x_i - y_i|.$$

6.2. *Bias bounds.* Since Huber [18] the median is known to be optimally robust against bias due to data contamination. Such robustness is essential to our analysis of the behavior of median filtering near edges. In effect, data contributed by the "other side" of an edge act as contamination that the median can optimally resist.

Consider now a composite dataset of $n + m$ points made from $x_i$, $i = 1, \ldots, n$, $y_i$, $i = 1, \ldots, m$. Think of the $y$'s as "bad" contamination of the "good" $x_i$ which may potentially corrupt the value of the median. How much damage can the $y$'s do? Equation (6.1) yields the mixture bounds

$$
\begin{aligned}
(6.3) \quad \mathrm{Med}_{n+m}(x_1, \ldots, x_n, -\infty, \ldots, -\infty) &\leq \mathrm{Med}_{n+m}(x_1, \ldots, x_n, y_1, \ldots, y_m) \\
&\leq \mathrm{Med}_{n+m}(x_1, \ldots, x_n, \infty, \ldots, \infty).
\end{aligned}
$$

Observe that if $m < n$, then the median of the combined sample cannot be larger than the maximum of the $x$'s nor can it be smaller than the minimum of the $x$'s:

$$
\mathrm{Med}_{n+m}(x_1, \ldots, x_n, -\infty, \ldots, -\infty) \geq \min(x_1, \ldots, x_n)
$$

and

$$
\mathrm{Med}_{n+m}(x_1, \ldots, x_n, \infty, \ldots, \infty) \leq \max(x_1, \ldots, x_n).
$$

Generalizing this observation leads to bias bounds employing the empirical quantiles of $x_1, \ldots, x_n$. Let $F_n(t) = n^{-1} \# \{i : x_i \leq t\}$ be the usual cumulative distribution function of the numbers $(x_i)$, and let $F_n^{-1}$ denote the empirical quantile function. Set $\varepsilon = m/(m + n)$ and suppose $\varepsilon \in (0, 1/2)$. As in [18], we bound the median of the combined sample by the quantiles of the "good" data only:

$$
(6.4) \quad F_n^{-1}\left( \frac{1/2 - \varepsilon}{1 - \varepsilon} \right) \leq \mathrm{Med}_{n+m}(x_1, \ldots, x_n, y_1, \ldots, y_m) \leq F_n^{-1}\left( \frac{1/2}{1 - \varepsilon} \right).
$$

This inequality will be helpful later, when the combined sample corresponds to all the data within a window of the median filter, the "good" data correspond to the part of the window on the "right" side of an edge, and the "bad" data correspond to the part of the window on the "wrong" side of the edge.

6.3. *Variance bounds for uncontaminated data.* The stochastic properties of the median are also crucial in our analysis; in particular we need bounds on the variance of the median of "uncontaminated" samples, that is, of the samples $(Z_i)_{i=1}^m$, $Z_i \overset{\text{i.i.d.}}{\sim} \Psi$. The following bounds on the variance of empirical medians behave similarly to expressions for variances of empirical averages.

LEMMA 6.1. *Suppose $\Psi$ satisfies* [Shape] *and* [Decay]. *Then, there are constants $C_1, C_2$ depending only on $\zeta(\Psi)$, such that*

$$\frac{C_1}{m} \leq \mathbb{E}[\operatorname{Med}_m(Z_1, \ldots, Z_m)^2] \leq \frac{C_2}{m}, \qquad m = 1, 2, \ldots.$$

PROOF. In Section 15, we prove Lemma 15.1 which states that [Shape] and [Decay] imply a condition due to David Mason, allowing us to apply Proposition 2 in [23]. □

We also need to analyze the properties of repeated medians (medians of medians). Borrowing ideas of Rousseeuw and Bassett [29], we prove the following in Section 15.

LEMMA 6.2. *Assume $\Psi$ satisfies* [Shape] *and* [Decay], *and consider $Z_1, \ldots, Z_m$ a sample from $\Psi$. Let $\Psi_m$ denote the distribution of $m^{1/2} \operatorname{Median}\{Z_1, \ldots, Z_m\}$. Then, for all $m$, $\Psi_m$ satisfies* [Shape] *and* [Decay]. *More precisely, there is a constant $C$ such that, for $m$ large enough, $\psi_m(x)(1 + |x|)^4 \leq C$ for all $x$.*

6.4. *Variance of empirical quantiles, uncontaminated data.* Because of the bias bound (6.4) it will be important to control not only the empirical median, but also other empirical quantiles besides $p = \frac{1}{2}$.

Let $Z_{m,p}$ denote the empirical $p$-quantile of $Z_1, \ldots, Z_m$, a sample from $\Psi$. That is, with $Z_{(1)}, \ldots, Z_{(m)}$ denoting order statistics of the sample, and $0 < p < 1$,

$$Z_{m,p} \equiv Z_{(1 + \lfloor mp \rfloor)}.$$

LEMMA 6.3. *Fix $\zeta > 1$. Let $\Psi$ satisfy* [Shape] *and* [Decay]. *Define*

$$(6.5) \qquad \alpha = \begin{cases} \dfrac{5\zeta - 3}{4\zeta - 4}, & \text{if } \zeta > 3, \\[2mm] \dfrac{\zeta}{\zeta - 1}, & \text{if } \zeta \leq 3. \end{cases}$$

*There is a constant $C > 0$ such that, for all sufficiently large positive integer $m$ and $p \in (2\alpha/m, 1 - 2\alpha/m)$,*

$$\mathbb{E}[Z_{m,p}^2] \leq C(p(1-p))^{-2\alpha+2}.$$

PROOF. Again noting Lemma 15.1, we are entitled to apply Proposition 2 in [23]. We then invoke Lemma 15.2. □

In words, provided that we do not consider quantiles $p$ very close to the extremes 0 and 1, the variance is well-controlled. The rate at which the variance blows up as $p \to 0$ or 1 is ultimately determined by the value of $\zeta > 1$ and will be of crucial significance for some bounds below.

6.5. *MSE lower bound for contaminated data.* As a final key ingredient in our analysis, we develop a simple lower bound on the mean-square displacement of the empirical median of contaminated data. Let

$$\hat{\mu}_{n,m,\Delta} \sim \mathrm{Med}_{n+m}(Z_1, \ldots, Z_n, Z_{n+1} + \Delta, \ldots, Z_{n+m} + \Delta).$$

In words this is the empirical median of $n + m$ values, the first $n$ of which are "good" data with median 0, and the last $m$ of which are contaminated, having median $\Delta$. For $\Delta > 0$ and $\varepsilon \in (0, 1)$, define the mixture CDF

$$(6.6) \qquad F_{\varepsilon,\Delta}(\cdot) = (1 - \varepsilon)\Psi(\cdot) + (1 - \varepsilon)\Psi(\cdot - \Delta).$$

Let $\mu = \mu(\varepsilon, \Delta)$ be the corresponding population median:

$$\mu = F_{\varepsilon,\Delta}^{-1}(\tfrac{1}{2}).$$

Actually, $\mu(\varepsilon, \Delta)$ is almost the population median of the empirical median $\hat{\mu}_{n,m,\Delta}$. More precisely, we have the following lemma.

LEMMA 6.4. *For $\varepsilon = m/(n + m)$,*

$$P\{\hat{\mu}_{n,m,\Delta} \geq \mu(\varepsilon, \Delta)\} \geq 1/2.$$

PROOF. For all $a \in \mathbb{R}$

$$\mathbb{P}\{\hat{\mu}_{n,m,\Delta} < a\} = \mathbb{P}\left\{\sum_{j=1}^{n}\{Z_j < a\} + \sum_{j=n+1}^{n+m}\{Z_j < a - \Delta\} \geq (n+m)/2\right\}.$$

Applying a result of Hoeffding [35], page 805, Inequality 1, we get that, for $a \leq \mu(\varepsilon, \Delta)$, the right-hand side is bounded by

$$\mathbb{P}\{\mathrm{Bin}(n + m, F_{\varepsilon,\Delta}(a)) \geq (n + m)/2\} \leq 1/2. \qquad \square$$

For each fixed $\Delta > 0$, increasing contamination only increases the population median:

$$(6.7) \qquad \mu(\varepsilon, \Delta) \text{ is an increasing function of } \varepsilon \in (0, 1/2).$$

Combining the last two observations, we have the MSE lower bound

$$(6.8) \qquad E\hat{\mu}_{n,m,\Delta}^2 \geq \mu(\varepsilon_0, \Delta)^2/2, \qquad \varepsilon_0 < m/(n + m).$$

**7. Proof of Theorem 2.1.** We now turn to proofs of our main results. In what follows, $C$ stands for a generic positive constant that depends only on the relevant function class and the distribution $\Psi$; its value may change from appearance to appearance. Also, to simplify the notation we use $\mathcal{W}(i) \equiv \mathcal{W}[n, h](i)$, $L_h(i) \equiv L_h[Y_n](i)$ and so on. We also write $\mathcal{F}_1$ in place of PLIP.

7.1. *Upper bound.* Fix $f \in \mathcal{F}_1$. Let $x_1, \ldots, x_N \in (0, 1)$ denote the points where $f$ is allowed to be discontinuous. Here and throughout the rest of the paper we assume that $h \geq 1/n$. Indeed, if to the contrary $h < 1/n$, then for all $i \in \mathbb{I}_n$, $\mathcal{W}(i) = \{i\}$ and so $\mathcal{R}_n(L_h; f) = \sigma^2 \asymp 1$.

We will demonstrate that

$$\mathcal{R}(L_h; \mathcal{F}_1) \leq C \left( h + \frac{1}{nh} \right). \tag{7.1}$$

Minimizing the right-hand side as a function of $h \geq 1/n$ gives $h_n = n^{-1/2}$, which implies our desired upper bound:

$$\mathcal{R}(L_h; \mathcal{F}_1) \leq C n^{-1/2}. \tag{7.2}$$

This upper bound may also be obtained from existing results, because $\mathcal{F}_1$ is included in a total-variation ball. Indeed,

$$\|f\|_{BV} \leq 2\beta + N,$$

so, in an obvious notation $\mathcal{F}_1 \subset BV(2\beta + N)$. From standard results on estimation of functions of bounded variation in white noise—[10, 21]—we know that

$$\inf_{h>0} \mathcal{R}(L_h; BV(\nu)) \leq \nu \cdot n^{-1/2},$$

which implies (7.2). Nevertheless, we spell out here an argument based on bias-variance trade-off, because this sort of trade-off will be used again repeatedly below.

Write the mean-squared error as squared bias plus variance: $\mathcal{R}_n(L_h; f) = B^2 + V$, where

$$B^2 = \frac{1}{n} \sum_{i=1}^{n} (\mathbb{E}[L_h(i)] - f(i/n))^2 \quad \text{and} \quad V = \frac{1}{n} \sum_{i=1}^{n} \operatorname{var}[L_h(i)].$$

For the variance, since the $\{Y_n(j) : j \in \mathbb{I}_n\}$ are pairwise uncorrelated and their variance is equal to $\sigma^2$, we have

$$\operatorname{var}[L_h(i)] = \frac{\sigma^2}{\#\mathcal{W}(i)} \leq \frac{\sigma^2}{nh}.$$

Therefore,

$$V \leq \frac{\sigma^2}{nh}. \tag{7.3}$$

For the bias, recall that $\mathbb{E}[Y_n(j)] = f(j/n)$ for all $j \in \mathbb{I}_n$, so

$$\mathbb{E}[L_h(i)] = \frac{1}{\#\mathcal{W}(i)} \sum_{j \in \mathcal{W}(i)} f(j/n).$$

Now consider separately cases where $i$ is "near to" and "far from" the discontinuity. Specifically, define:

- Near: $\Delta = \{i \in \mathbb{I}_n : \min_s |i/n - x_s| \leq h\}$—in words, $\Delta$ is the set of points where the smoothing window *does* meet the discontinuity.
- Far: $\Delta^c = \{i \in \mathbb{I}_n : \min_s |i/n - x_s| > h\}$—$\Delta^c$ is the set of points where the smoothing window *does not* meet the discontinuity.

Near the discontinuity, use the fact that as $f$ takes values in $[0,1]$,

$$(7.4) \qquad |\mathbb{E}[L_h(i)] - f(i/n)| \leq 1.$$

Far from the discontinuity, we apply a sharper estimate that we now develop. Let $i \in \Delta^c$ and consider $j \in \mathcal{W}(i)$. The local Lipschitz constant bound $\beta$ gives

$$(7.5) \qquad |f(j/n) - f(i/n)| \leq h \sup_{x \in [i/n, j/n]} L_x(f) \leq \beta h,$$

which implies

$$(7.6) \qquad |\mathbb{E}[L_h(i)] - f(i/n)| \leq \beta h, \qquad i \in \Delta^c.$$

Combining (7.4) and (7.6), we bound the squared bias by

$$(7.7) \qquad B^2 \leq \frac{\#\Delta^c}{n} \beta^2 h^2 + \frac{\#\Delta}{n}.$$

The number of "near" terms obeys $0 \leq \#\Delta \leq N(2nh + 1)$, and of course the fraction of "far" terms obeys $\#\Delta^c/n \leq 1$, so we get

$$(7.8) \qquad B^2 \leq \beta^2 h^2 + 2Nh + N/n \leq Ch,$$

since $h \geq 1/n$; we may take $C = (\beta^2 + 3N)$.

Hence, $\mathcal{R}(L_h; f) \leq C(h + 1/(nh))$, and this bound does not depend on $f \in \mathcal{F}_1$, so (7.1) follows.

7.2. *Lower bound.* Let $f$ be the indicator function of the interval $[1/2, 1]$. Then $f \in \mathcal{F}_1$ for all $N \geq 1$ and $\beta > 0$.

For the variance, since $\#\mathcal{W}(i) \leq 3nh$, we have

$$V \geq \frac{\sigma^2}{3nh}.$$

For the squared bias, we show that the pointwise bias is large near the discontinuity. For example, take $n/2 - nh/2 \leq i < n/2$, so that $f(i/n) = 0$ and therefore

$$|\mathbb{E}[L_h(i)] - f(i/n)| = \mathbb{E}[L_h(i)] = \frac{\#\{j \in \mathcal{W}(i) : j/n \geq 1/2\}}{\#\mathcal{W}(i)}.$$

Since $\#\{j \in \mathcal{W}(i) : j/n \geq 1/2\} \geq nh/2$ and $\#\mathcal{W}(i) \leq 3nh$, the pointwise bias exceeds $1/6$:

$$|\mathbb{E}[L_h(i)] - f(i/n)| \geq \frac{nh/2}{3nh} \geq 1/6.$$

Therefore,

$$B^2 \geq \frac{\#\{i \in \mathbb{I}_n : n/2 - nh/2 \leq i < n/2\}}{n}(1/6)^2 \geq Ch,$$

where $C = 1/72$ will do.

Combining bias and variance bounds, we have for any choice of radius $h$,

$$\mathcal{R}_n(L_h; f) \geq C\left(h + \frac{1}{nh}\right) \geq Cn^{-1/2}.$$

**8. Proof of Theorem 2.2.**  Now we analyze median filtering. The proof parallels that for Theorem 2.1 and uses the results of Section 6. Here too, we use abbreviated notation.

8.1. *Upper bound.*  Fix $f \in \mathcal{F}_1$. Let $x_1, \ldots, x_N \in (0,1)$ be the points where $f$ may be discontinuous. Without loss of generality, we again let $h \geq 1/n$.

We will show that

$$(8.1) \qquad \mathcal{R}_n(M_h; \mathcal{F}_1) \leq C\left(h + \frac{1}{nh}\right).$$

As in the proof of Theorem 2.1, picking $h_n = n^{-1/2}$ in (8.1) implies our desired upper bound, namely:

$$\inf_{h>0} \mathcal{R}_n(M_h; \mathcal{F}_1) \leq Cn^{-1/2}.$$

To get started, we invoke the monotonicity and Lipschitz properties of the median (6.1)–(6.2), yielding

$$|M_h(i) - f(i/n)| \leq \max_{j \in \mathcal{W}(i)} |f(j/n) - f(i/n)| + \sigma|\hat{Z}(i)|,$$

where $\hat{Z}(i) = \text{Median}\{Z_n(j) : j \in \mathcal{W}(i)\}$.

Near the discontinuity, we again observe that since $f$ takes values in $[0,1]$, we have

$$\max_{j \in \mathcal{W}(i)} |f(j/n) - f(i/n)| \leq 1$$

for all $i \in \mathbb{I}_n$, and so

$$(8.2) \qquad |M_h(i) - f(i/n)| \leq 1 + \sigma|\hat{Z}(i)| \qquad \forall i \in \mathbb{I}_n.$$

Now consider the set $\Delta^c = \{i \in \mathbb{I}_n : \min_s |i/n - x_s| > h\}$ far from the discontinuity. Using (7.5), we get

$$(8.3) \qquad |M_h(i) - f(i/n)| \leq \beta h + \sigma|\hat{Z}(i)| \qquad \forall i \in \Delta^c.$$

Using (8.2) and (8.3), we get

$$(8.4) \qquad \mathcal{R}_n(M_h; f) \leq \frac{\#\Delta^c}{n}\beta^2 h^2 + \frac{\#\Delta}{n} + \frac{1}{n}\sum_{i=1}^{n} \mathbb{E}[\hat{Z}(i)^2].$$

The term on the far right is a variance term, which can be handled using Lemma 6.1 at sample size $m = \#\mathcal{W}(i) \geq nh$, yielding variance $V \leq C/(nh)$. The bias terms involve $\#\Delta$ and $\#\Delta^c$ and are completely analogous to the case of linear filtering and are handled just as at (7.8), using $0 \leq \#\Delta \leq N(2nh + 1)$. We obtain $\mathcal{R}_n(M_h; f) \leq C(h + 1/(nh))$. Since this does not depend on $f \in \mathcal{F}_1$, (8.1) follows.

8.2. *Lower bound.* Let $f$ be the indicator function of the interval $[1/2, 1]$. Surely, $f \in \mathcal{F}_1$ for any $N \geq 1$ and $\beta > 0$.

For $i \in \Delta^c$,

$$|M_h(i) - f(i/n)| = \sigma|\hat{Z}(i)|,$$

so that, by Lemma 6.1,

(8.5)  $$\mathbb{E}[(M_h(i) - f(i/n))^2] \geq C\frac{1}{nh}.$$

Therefore,

$$\frac{1}{n}\sum_{i \in \Delta^c}\mathbb{E}[(M_h(i) - f(i/n))^2] \geq C\frac{1}{nh}.$$

For $i \in \Delta$, we view the window as consisting of a mixture of "good" data, on the same side of the discontinuity as $i$ together with "bad" data, on the other side. Thus

$$|M_h(i) - f(i/n)| = \sigma|K_n(i)|,$$

where, with $w(i) = \#\mathcal{W}(i) \asymp nh$, and

$$\rho(i) = \#\{j : j \in \mathcal{W}(i) \text{ and on the same side of the discontinuity}\}$$

we have

$$K_n(i) \sim \text{Median}\{Z_1, \ldots, Z_{\rho(i)w(i)}, Z_{1+\rho(i)w(i)} + 1/\sigma, \ldots, Z_{w(i)} + 1/\sigma\}.$$

This is exactly a median of contaminated data as discussed in Section 6.5, with $\Delta = 1/\sigma$, $m + n = w(i)$, $n = \rho(i)w(i)$ and $\varepsilon = 1 - \rho(i)$. Invoking Lemma 6.4 with $\mu_i = \mu(1 - \rho(i), 1/\sigma)$, applying (6.8) for $\varepsilon_0 = 1/5$ and setting $C = \sigma^2\mu^2(1/5, 1/\sigma)/2$, we have

(8.6)  $$\mathbb{E}[(M_h(i) - f(i/n))^2] \geq C \qquad \forall i \text{ such that } \rho(i) \leq 4/5.$$

Since $\#\{i : \rho(i) \leq 4/5\} \asymp nh$, we get

$$\frac{1}{n}\sum_{i \in \Delta}\mathbb{E}[(M_h(i) - f(i/n))^2] \geq Ch.$$

Combining pieces, we get for any choice of $h$

$$\mathcal{R}_n(M_h; f) \geq C\left(h + \frac{1}{nh}\right) \geq C \cdot n^{-1/2},$$

which matches the upper bound.

**9. Proof of Theorem 4.1.** We turn to the setting of Section 4: asymptotically negligible noise per pixel: $\sigma = \sigma_n = o(1)$.

9.1. *Linear filtering.* By just carrying the variance term in Section 7 and comparing with the no-smoothing rate, we immediately see that

$$\mathcal{R}(L_h; \text{SEP-PLIP}) \asymp \sigma_n n^{-1/2} \wedge \sigma_n^2.$$

Hence, linear filtering improves on no-smoothing if, and only if, $\sigma_n n^{1/2} \to \infty$.

9.2. *Upper bound for median filtering.* We refine the argument from Section 8 in the case where the discontinuities are well-separated. Let $\mathcal{F}_1^+ = \text{SEP-PLIP}$ as defined in Section 5. Note that (9.7) will be used in the proof of Theorem 5.1.

Assuming $n > 2/\eta$, choose $h \in [1/n, \eta/2)$. Since median filtering is local and the discontinuities are $\eta$-separated, we may assume that $f$ only has $N = 1$ discontinuity point $x_1$.

Far from the discontinuity, at $i \in \Delta^c$, we use (8.2) and Lemma 6.1 to get

$$(9.1) \qquad \mathbb{E}[(M_h(i) - f(i/n))^2] \leq C\left(h^2 + \frac{\sigma_n^2}{nh}\right), \qquad i \in \Delta^c;$$

note that $\sigma_n$ is now nonconstant.

Near the discontinuity, we now take more seriously the viewpoint that the window contains "good" data (on the same side of the discontinuity) and "bad" data (on the other side) and we bound the MSE more carefully than before.

So take $i \in \Delta$. Define the "good" subset $\mathcal{G}(i)$—the subset of the window on the same side of the discontinuity as $i$—by

$$(9.2) \qquad \mathcal{G}(i) = \begin{cases} \mathcal{W}(i) \cap [nx_1, n], & i/n > x_1, \\ \mathcal{W}(i) \cap [1, nx_1), & i/n < x_1. \end{cases}$$

Let $\rho(i) = \#\mathcal{G}(i)/\#\mathcal{W}(i)$ and $\varepsilon(i) = 1 - \rho(i)$. The window $\mathcal{W}(i)$ provides an $\varepsilon$-contaminated sample in the sense of Huber [18].

By the contamination bias bound (6.4), $M_h(i)$ lies between the $\frac{1/2-\varepsilon}{1-\varepsilon} = (2\rho(i)-1)/(2\rho(i))$ and $\frac{1}{2(1-\varepsilon)} = 1/(2\rho(i))$ quantiles of $\{Y_n(j) : j \in \mathcal{G}(i)\}$. Also, as in (7.5), we have

$$|f(j/n) - f(i/n)| \leq \beta h \qquad \forall j \in \mathcal{G}(i),$$

so the quantiles of these "good data" $\{Y_n(j) : j \in \mathcal{G}(i)\}$ are small perturbations of the quantiles of corresponding zero-median data $\{Z_n(j) : j \in \mathcal{G}(i)\}$. Let $Q_n(i)$ denote the maximum absolute value of the empirical $(2\rho(i) - 1)/(2\rho(i))$ and $1/(2\rho(i))$ quantiles of $\{Z_n(j) : j \in \mathcal{G}(i)\}$. By (6.2) and (6.4)

$$(9.3) \qquad |M_h(i) - f(i/n)| \leq \beta h + \sigma_n Q_n(i).$$

We combine (9.3) with (8.2) and Lemma 6.1 to get

$$\mathbb{E}[(M_h(i) - f(i/n))^2] \leq C\mathbb{E}[(1 + \sigma_n^2 \hat{Z}(i)^2) \wedge (\beta^2 h^2 + \sigma_n^2 Q_n(i)^2)]$$
$$\leq C(1 + \sigma_n^2 \mathbb{E}[\hat{Z}(i)^2]) \wedge (\beta^2 h^2 + \sigma_n^2 \mathbb{E}[Q_n(i)^2])$$
$$\leq C(h^2 + 1 \wedge \sigma_n^2 \mathbb{E}[Q_n(i)^2]),$$

where for the last inequality we used Lemma 6.1 together with $h \geq 1/n$, and the fact that $a \wedge (b+c) \leq (a \wedge b) + c$ for any $a, b, c \geq 0$.

Recalling Section 6.4, let $Z_{m,p}$ denote the empirical $p$-quantile of $Z_1, \ldots, Z_m$, a sample from $\Psi$. Because $\Psi$ is symmetric about 0, $|Q_n(i)|$ is stochastically majorized by $2|Z_{m(i),p(i)}|$, with $m(i) = \#\mathcal{G}(i) = \rho(i)\#\mathcal{W}(i)$ and $p(i) = 1/(2\rho(i))$. Hence

$$(9.4) \qquad \mathbb{E}[Q_n(i)^2] \leq 4\mathbb{E}[Z_{m(i),p(i)}^2].$$

Using (9.4) and Lemma 6.3, we obtain, for $i \in \Delta$,

$$\mathbb{E}[(M_h(i) - f(i/n))^2] \leq C(h^2 + 1 \wedge \sigma_n^2(p(i)(1-p(i)))^{-2\alpha+2}).$$

Therefore,

$$(9.5) \quad \mathcal{R}_n(M_h; f) \leq C\left(h^2 + \frac{\sigma_n^2}{nh} + \frac{1}{n}\sum_{i \in \Delta} 1 \wedge \sigma_n^2(p(i)(1-p(i)))^{-2\alpha+2}\right).$$

We focus on the last term on the right-hand side.

Let $\delta(i) = |i/n - x_1|$; since $i \in \Delta$, $\delta(i) \leq h$. We have

$$\rho(i) = \frac{[n\delta(i)] + [nh] + 1}{2[nh] + 1} \geq \frac{1}{2} + C\frac{\delta(i)}{h}.$$

So there is a constant $C > 0$ such that

$$(9.6) \qquad\qquad p(i) \leq 1 - C\delta(i)/h \qquad \forall i \in \Delta.$$

Note that we always have $p(i) \geq 1/2$.

Therefore

$$\frac{1}{n}\sum_{i \in \Delta} 1 \wedge \sigma_n^2(p(i)(1-p(i)))^{-2\alpha+2} \leq C\frac{1}{n}\sum_{i \in \Delta} 1 \wedge \sigma_n^2(\delta(i)/h)^{-2\alpha+2}$$

$$\leq Ch\frac{1}{nh}\sum_{i=1}^{nh} 1 \wedge \sigma_n^2(i/(nh))^{-2\alpha+2}$$

$$\leq Ch\int_0^1 1 \wedge \sigma_n^2 s^{-2\alpha+2}\, ds = Ch\nu_n$$

with

$$\nu_n = \begin{cases} \sigma_n^2, & \text{if } \zeta > 3, \\ \sigma_n^2 \log(1/\sigma_n), & \text{if } \zeta = 3, \\ \sigma_n^{\zeta-1}, & \text{if } \zeta < 3. \end{cases}$$

Combining pieces gives

$$(9.7) \qquad \mathcal{R}_n(M_h; f) \leq C\left(h^2 + \frac{\sigma_n^2}{nh} + h\nu_n\right).$$

Optimizing the right-hand side over $h$, we get

$$\mathcal{R}_n(M_h; f) \leq C(\nu_n^{1/2} \vee \sigma_n^{1/3}n^{-1/6}) \cdot \sigma_n n^{-1/2} = o(\sigma_n n^{-1/2}).$$

This bound improves on no-smoothing if $\sigma_n n \to \infty$.

9.3. *Lower bound for median filtering.* A lower bound is not needed to prove Theorem 4.1. However, we will use the following lower bound in the proof of Theorem 5.1.

Let $f$ be the indicator function of the interval $[1/2, 1]$. A lower bound is obtained by using the arguments in Section 8.2, but this time carrying $\sigma_n$ along and noticing that $\mu(\varepsilon, \Delta)$ is increasing in $\Delta$. One gets

$$(9.8) \qquad \mathcal{R}_n(M_h; f) \geq C\left(h\sigma_n^2 + \frac{\sigma_n^2}{nh}\right).$$

This bound matches the upper bound, for example, when $\zeta > 3$ and $\sigma_n n^{1/4} \to \infty$. This is the setting that will arise in Section 13.

## 10. Proof of Theorem 3.1.

We consider two-dimensional linear filtering. The structure of the argument parallels the one-dimensional case presented in Section 7. The main difference involves counting points near to discontinuities.

10.1. *Upper bound.* We also write $\mathcal{F}_2$ in place of CPLIP. Fix $f \in \mathcal{F}_2$. We call $\gamma_1, \ldots, \gamma_N \in \Gamma$ the curves where $f$ may be discontinuous.

As before, we write $MSE = B^2 + V$.

Again, we may assume $h \geq 1/n$. Since $\#\mathcal{W}(\mathbf{i}) \geq (nh)^2$, we have $V \leq 1/(nh)^2$.

In the two-dimensional case, we define proximity to singularity as follows. Write $d(A, B)$ for Haussdorff distance between subsets $A$ and $B$ of the unit square:

- **Far:** Let $\Delta^c = \{\mathbf{i} \in \mathbb{I}_n^2 : \min_s d(\mathbf{i}/n, \gamma_s) > h\}$.
- **Near:** Let $\Delta = \{\mathbf{i} \in \mathbb{I}_n^2 : \min_s d(\mathbf{i}/n, \gamma_s) \leq h\}$.

Using the exact same arguments as in Section 7, we obtain the equivalent of (7.7):

$$B^2 \leq \frac{\#\Delta^c}{n^2}\beta^2 h^2 + \frac{\#\Delta}{n^2} \leq \beta^2 h^2 + \frac{\#\Delta}{n^2}.$$

Lemma 16.1 provides an estimate for $\#\Delta$ which, when used in the above expression, implies $B^2 \leq Ch$.

Thus, $\mathcal{R}_n(L_h; f) \leq C(h + 1/(nh)^2)$. The right-hand side does not depend on $f \in \mathcal{F}_2$, so

$$\mathcal{R}_n(L_h; \mathcal{F}_2) \leq C\left(h + \frac{1}{(nh)^2}\right).$$

Minimizing the right-hand side over $h \geq 1/n$ gives $h = n^{-2/3}$, yielding $\mathcal{R}(L_h; \mathcal{F}_2) \leq Cn^{-2/3}$.

10.2. *Lower bound.* Fix $0 < \zeta < (1/2) \wedge (\lambda/4)$ and let $f$ be the indicator function of the axis-aligned square of sidelength $\zeta$ centered at $(1/2, 1/2)$, namely $f = 1_S$ where $S = [1/2 - \zeta/2, 1/2 + \zeta/2] \times [1/2 - \zeta/2, 1/2 + \zeta/2]$. Certainly, $f \in \mathcal{F}_2$.

Again, the variance $V \geq 1/(3nh)^2$. For the squared bias $B^2$, we show that the pointwise bias is of order 1 near the discontinuity. For example, take $\mathbf{i} \in \mathbb{I}_n^2$ such that

$$\mathbf{i}/n \in [1/2 - \zeta/2 - h/2, 1/2 - \zeta/2] \times [1/2 - \zeta/2, 1/2 + \zeta/2],$$

so that $f(\mathbf{i}/n) = 0$ and therefore the bias obeys

$$|\mathbb{E}[L_h(\mathbf{i})] - f(\mathbf{i}/n)| = \mathbb{E}[L_h(\mathbf{i})] = \frac{\#\{\mathbf{j} \in \mathcal{W}(\mathbf{i}) : \mathbf{j} \in nS\}}{\#\mathcal{W}(\mathbf{i})}.$$

For such $\mathbf{i}$, $\#\{\mathcal{W}(\mathbf{i}) \cap nS\}$ is of order $(nh)^2$, since the intersection of the disc of radius $h$ centered at $\mathbf{i}/n$ with $S$ contains a square of sidelength $Ch$. Therefore, the bias is of order 1 for such $\mathbf{i}$, and there are order $n^2h$ such $i$. Hence the squared bias $B^2$ is at least of order $h$.

Combining pieces, we get for all $h$,

$$\mathcal{R}_n(L_h; f) \geq C \cdot \left(h + \frac{1}{(nh)^2}\right) \geq C \cdot n^{-2/3}.$$

**11. Proof of Theorem 3.2.** The structure of the proof is identical to the case of one-dimensional signals presented in Section 8. In the details, the only significant difference is on computing the number of points away from discontinuities. We use the same definitions $\Delta$ and $\Delta^c$ as in Section 8.

11.1. *Upper bound.* Fix $f \in \mathcal{F}_2$. We call $\gamma_1, \ldots, \gamma_N \in \Gamma$ the curves where $f$ may be discontinuous. Again, we may assume $h \geq 1/n$.

Define $\Delta^c = \{\mathbf{i} \in \mathbb{I}_n^2 : \min_s d(\mathbf{i}/n, \gamma_s) > h\}$. Using the exact same arguments as in Section 8, we obtain the equivalent of (8.4):

$$\mathcal{R}_n(M_h; f) \leq \frac{\#\Delta^c}{n^2}\beta^2 h^2 + \frac{\#\Delta}{n^2} + \frac{1}{n^2}\sum_{\mathbf{i} \in \mathbb{I}_n^2}\mathbb{E}[\hat{Z}(\mathbf{i})^2].$$

Lemma 16.1 provides an estimate for $\#\Delta$ which, when used in the above expression, leads to $\mathcal{R}_n(L_h; f) \leq C(h + 1/(nh)^2)$. From there we conclude as in Section 10.1.

### 11.2. *Lower bound.*

Let $f \in \mathcal{F}_2$ be the indicator function of a disc $D$. For $\mathbf{i} \in \Delta^c$, the equivalent of (8.5) holds and together with Lemma 16.1 implies

$$\frac{1}{n^2} \sum_{\mathbf{i} \in \Delta^c} \mathbb{E}[(M_h(\mathbf{i}) - f(\mathbf{i}/n))^2] \geq C \frac{1}{(nh)^2}.$$

The equivalent of (8.6) holds as well and implies

$$\frac{1}{n^2} \sum_{\mathbf{i} \in \Delta} \mathbb{E}[(M_h(\mathbf{i}) - f(\mathbf{i}/n))^2] \geq C \frac{1}{n^2} \#\{\mathbf{i} : \rho(\mathbf{i}) \leq 4/5\}.$$

We now show that, for $\mathbf{i}/n \in D^c$ such that $\delta(\mathbf{i}) \geq 1/n$, $\rho(\mathbf{i}) \leq 1/2 + C\delta(\mathbf{i})/h$. Let $y \in \partial D$ be the closest point to $\mathbf{i}/n$ and $L$ the tangent to $\partial D$ at $y$. $L$ divides $B(\mathbf{i}/n, h)$ into two parts $A$ and $B(\mathbf{i}/n, h) \cap A^c$, where $\mathbf{i}/n \in A$. We have $A = A_0 \cup H$, where $H$ is the open half disc with diameter parallel to $L$ that does not intersect $\partial D$. We have $\mathcal{G}(\mathbf{i}) = \mathbb{I}_n^2 \cap nA$, so that $\#\mathcal{G}(\mathbf{i}) \leq \#\mathcal{W}(\mathbf{i})/2 + \#(\mathbb{I}_n^2 \cap nA_0)$. $A_0$ is contained within a rectangular region $R$ with dimensions $\delta(\mathbf{i})$ by $2h$ and for any rectangular region, $|R| \leq C|R|n^2 + O(nh)$. Hence, since $n\delta(\mathbf{i}) \geq 1$,

$$\#(\mathbb{I}_n^2 \cap nA_0) \leq \#(\mathbb{I}_n^2 \cap nR) \leq Ch\delta(\mathbf{i})n^2.$$

Therefore, $\rho(\mathbf{i}) \leq 1/2 + C\delta(\mathbf{i})/h$.

We thus have

$$\#\{\mathbf{i} : \rho(\mathbf{i}) \leq 4/5\} \geq \#\{\mathbf{i} : 1/n \leq \delta(\mathbf{i}) \leq Ch\}.$$

Let $K = |\{x : d(x, \partial D) \leq Ch\}|$. We have $nh \gg 1$ so that

$$K \subset \bigcup_{\mathbf{i}/n \in R} B(\mathbf{i}/n, 2/n),$$

which implies $|K| \leq C\#\{\mathbf{i} : \delta(\mathbf{i}) \leq Ch\}/n^2$. By elementary calculus, $|K| \geq Ch$, so

$$\#\{\mathbf{i} : \rho(\mathbf{i}) \leq 4/5\} \geq Cn^2 h.$$

We obtain for all $h$

$$\mathcal{R}_n(M_h; f) \geq C\left(\frac{1}{(nh)^2} + h\right) \geq Cn^{-2/3}.$$

### 12. Proof of Theorem 4.2.

We consider again median filtering in the negligible-noise-per-pixel case of Section 4, this time in the two-dimensional setting.

12.1. *Linear filtering.* By just carrying the variance term in Section 10 and comparing with the no-smoothing rate, we immediately see that

$$\mathcal{R}(L_h; \text{SEP-CPLIP}) \asymp \sigma_n^{2/3} n^{-2/3} \wedge \sigma_n^2.$$

Hence, linear filtering improves on no-smoothing if, and only if, $\sigma_n n^{1/2} \to \infty$.

12.2. *Upper bound for median filtering.* We refine our arguments in the setting of asymptotically negligible noise level $\sigma = \sigma_n = o(1)$ with well-separated discontinuities. Let $\mathcal{F}_2^+ = \text{SEP-CPLIP}$ as defined in Section 5. Note that (12.1) will be used in the proof of Theorem 5.2.

Letting $n > 2/\eta$, choose $h \in [1/n, \eta/2)$. Since median filtering is local and the discontinuities are at least $\eta$ apart, we may assume that $N = 1$, namely that $f$ only has one discontinuity curve $\gamma \in \Gamma^+$. Being a Jordan curve, $\gamma$ partitions $[0,1]^2$ into two regions, the inside ($\Omega$) and the outside ($\Omega^c$).

We proceed as in Section 8, introducing $\delta(\mathbf{i}) = d(\mathbf{i}/n, \gamma)$ and

$$\mathcal{G}(\mathbf{i}) = \begin{cases} \mathcal{W}(\mathbf{i}) \cap \Omega, & \text{if } \mathbf{i}/n \in \Omega, \\ \mathcal{W}(\mathbf{i}) \cap \Omega^c, & \text{if } \mathbf{i}/n \in \Omega^c, \end{cases}$$

together with $\rho(\mathbf{i}) = \#\mathcal{G}(\mathbf{i})/\#\mathcal{W}(\mathbf{i})$ and $p(\mathbf{i}) = 1/(2\rho(\mathbf{i}))$.

Using the exact same arguments as in Section 8, we obtain the equivalent of (9.5):

$$\mathcal{R}_n(M_h; f) \leq C\left(h^2 + \frac{\sigma_n^2}{(nh)^2} + \frac{1}{n^2}\sum_{\mathbf{i}\in\Delta} 1 \wedge \sigma_n^2(p(\mathbf{i})(1-p(\mathbf{i})))^{-2\alpha+2}\right).$$

We bound the last term on the right-hand side by

$$\frac{1}{n^2}\sum_{\mathbf{i}\in\Delta\cap\Delta_1^c} 1 \wedge \sigma_n^2(p(\mathbf{i})(1-p(\mathbf{i})))^{-2\alpha+2} + \frac{1}{n^2}\sum_{\mathbf{i}\in\Delta_1} 1,$$

where $\Delta_1^c = \{\mathbf{i}\in\mathbb{I}_n^2 : \delta(\mathbf{i}) > 2(C_1 h^2 + n^{-1})\}$, the constant $C_1 > 0$ being given by Lemma 16.2—the second term in this last expression represents the bias due to the curvature of the discontinuity.

For $\ell = 0, \ldots, [nh]$, define

$$\Xi_\ell = \{\mathbf{i}\in\mathbb{I}_n^2 : \ell \leq n\delta(\mathbf{i}) < \ell+1\}.$$

We use Lemma 16.4 to get

$$\frac{1}{n^2}\sum_{\mathbf{i}\in\Delta_1} 1 = \frac{1}{n^2}\sum_{\ell=0}^{2n(C_1 h^2 + n^{-1})} \#\Xi_\ell \leq C\frac{1}{n}\sum_{\ell=0}^{2n(C_1 h^2 + n^{-1})} 1 \leq C(h^2 \vee n^{-1}).$$

We use Lemmas 16.2 and 16.4, and replicate the computations below (9.6) to get

$$\frac{1}{n^2} \sum_{\mathbf{i} \in \Delta \cap \Delta_1^c} 1 \wedge \sigma_n^2 (p(\mathbf{i})(1 - p(\mathbf{i})))^{-2\alpha+2} \leq C \frac{1}{n^2} \sum_{\ell=0}^{nh} \# \Xi_\ell \cdot (1 \wedge \sigma_n^2 (\ell/(nh))^{-2\alpha+2})$$

$$\leq C \frac{1}{n} \sum_{\ell=0}^{nh} 1 \cdot (1 \wedge \sigma_n^2 (\ell/(nh))^{-2\alpha+2})$$

$$\leq C h \nu_n.$$

Combining inequalities,

$$(12.1) \qquad \mathcal{R}_n(M_h; f) \leq C\left(h^2 + \frac{\sigma_n^2}{(nh)^2} + h\nu_n\right).$$

Optimizing the right-hand side over $h$, we get

$$\mathcal{R}_n(M_h; f) \leq C(\nu_n^{2/3} \vee \sigma_n^{1/3} n^{-1/3}) \cdot \sigma_n^{2/3} n^{-2/3} = o(\sigma_n^{2/3} n^{-2/3}).$$

This bound improves on no-smoothing if $\sigma_n n \to \infty$.

12.3. *Lower bound for median filtering.* This is not needed to prove Theorem 4.2. However, we will use the following lower bound in the proof of Theorem 5.2.

Let $f$ be the indicator function of a disc $D$ such that $f \in \mathcal{F}_2^+$. The low-noise-per-pixel case comes again from carrying $\sigma_n$ along, yielding

$$(12.2) \qquad \mathcal{R}_n(M_h; f) \geq C\left(h\sigma_n^2 + \frac{\sigma_n^2}{(nh)^2}\right).$$

This bound matches the upper bound, for example, when $\zeta > 3$ and $\sigma_n n^{1/3} \to \infty$. This is the setting that will arise in Section 14.

## 13. Proof of Theorem 5.1.

13.1. *Upper bound.* Without loss of generality, fix $\sigma = 1$. Let $1/n \leq h_1 < h_2 < 1$ to be chosen later as functions of $n$. We only need consider $h_1 \gg 1/n$, for otherwise the first pass does not reduce the noise level significantly. Also, for simplicity we assume that both $n_1 = h_1^{-1}$ and $nh_1$ are integers.

Fix $f \in \mathcal{F}_1^+$. Again, we may assume that $N = 1$ without loss of generality. Call $x_1 \in (0, 1)$ the point where $f$ is discontinuous and let

$$\{k_1\} = \{k : x_1 \in (B_k^{h_1}/n)\}.$$

Using (6.1)–(6.2), we have

$$Y_n^{h_1}(k) = f(kh_1) + Z_n^{h_1}(k) + V_n^{h_1}(k),$$

where

$$Z_n^{h_1}(k) = \text{Median}\{Z_n(i) : i \in B_k^{h_1}\},$$

$$|V_n^{h_1}(k)| \le U_n^{h_1}(k) = \max_{i \in B_k^{h_1}} |f(i/n) - f(kh_1)|.$$

Since for $k \ne k_1$, $f$ is locally Lipschitz in $B_k^{h_1}$, we have

$$|U_n^{h_1}(k)| \le \begin{cases} \beta h_1, & k \ne k_1, \\ 1, & k = k_1. \end{cases}$$

Let $\sigma_n' = (2nh_1 + 1)^{-1/2}$ and define $Z_n'(k) = Z_n^{h_1}(k)/\sigma_n'$. The $Z_n'(k)$'s are independent and identically distributed, and by Lemma 6.2, their distribution $\Psi_n'$ satisfies [Shape] and [Decay] with $\zeta(\Psi_n') \ge 4$ and implicit constant independent of $n$. Define $Y_n'(k) = f(kh_1) + \sigma_n' Z_n'(k)$. For $i \in B_k^{h_1}$, we have

$$\begin{aligned} |M^{h_1,h_2}[Y_n](i) &- f(i/n)| \\ &\le |M_{h_2}[Y_n^{h_1}](k) - M_{h_2}[Y_n'](k)| \\ &\quad + |M_{h_2}[Y_n'](k) - f(kh_1)| + |f(kh_1) - f(i/n)| \\ &\le |M_{h_2}[Y_n'](k) - f(kh_1)| + 2|U_n^{h_1}(k)|. \end{aligned}$$

Hence, using the bounds on $U_n^{h_1}(k)$, we have

$$\begin{aligned} \mathcal{R}_n(M^{h_1,h_2}; f) &= \frac{1}{n} \sum_{k \in \mathbb{I}_{n_1}} \sum_{i \in B_k^{h_1}} \mathbb{E}[(M^{h_1,h_2}[Y_n](i) - f(i/n))^2] \\ &\le C \frac{1}{n_1} \sum_{k \in \mathbb{I}_{n_1}} \mathbb{E}[(M^{h_2}[Y_n'](k) - f(k/n_1))^2] + U_n^{h_1}(k)^2 \\ &\le C \frac{1}{n_1} \sum_{k \in \mathbb{I}_{n_1}} \mathbb{E}[(M^{h_2}[Y_n'](k) - f(k/n_1))^2] + Ch_1. \end{aligned}$$

Using the upper bound (9.7) on the first term, which we may use since we are back to the original situation, we get

$$\mathcal{R}_n(M^{h_1,h_2}; f) \le C\left(h_2^2 + \frac{(\sigma_n')^2}{n_1 h_2} + h_2(\sigma_n')^2\right) + Ch_1.$$

We then replace $\sigma_n'$ by its definition $(2nh_1 + 1)^{-1/2}$ and minimize over $h_1$ and $h_2$, with $h_1 = n^{-2/3}$ and $h_2 = n^{-1/3}$, and obtain the desired upper bound valid for any $f \in \mathcal{F}_1^+$.

13.2. *Lower bound.* Fix $h_1 < h_2$. We again assume for convenience that $n_1 = 1/h_1$ and $nh_1$ are integers. Let $f$ be the indicator function of the interval $[t, 1]$, where $t$ is the middle of the unique interval of the form $[kh_1, (k+1)h_1]$ containing $1/2$. By the definition of $f$,

$$\frac{\#\{i \in B_{k_1}^{h_1} : f(i/n) = 0\}}{\# B_{k_1}^{h_1}} \in [1/3, 2/3].$$

Hence, because $M^{h_1,h_2}[Y_n](i) = M^{h_1,h_2}[Y_n](j)$ for all $i, j \in B_{k_1}^{h_1}$,

$$\frac{1}{n} \sum_{i \in B_{k_1}^{h_1}} \mathbb{E}[(M^{h_1,h_2}[Y_n](i) - f(i/n))^2] \geq C \frac{nh_1}{n} = Ch_1.$$

Now, because $U_n^{h_1}(k) = 0$ if $k \neq k_1$, we have

$$\frac{1}{n} \sum_{i \notin B_{k_1}^{h_1}} \mathbb{E}[(M^{h_1,h_2}[Y_n](i) - f(i/n))^2] = \frac{1}{n_1} \sum_{k \neq k_1} \mathbb{E}[(M_{h_2}[Y_n'](k) - f(k/n_1))^2].$$

We then use (9.8), which applies the same here even though we omit $k = k_1$.

Combining the cases $k = k_1$ and $k \neq k_1$, we get

$$\mathcal{R}_n(M^{h_1,h_2}; f) \geq C \left( \frac{(\sigma_n')^2}{n_1 h_2} + h_2 (\sigma_n')^2 \right) + Ch_1.$$

We conclude by noticing that the right-hand side is larger than $n^{-2/3}$ for all choices of $h_1 < h_2$.

## 14. Proof of Theorem 5.2.

14.1. *Upper bound.* We follow the line of arguments in Section 13.

Fix $f \in \mathcal{F}_2^+$. Again, we may assume that $N = 1$ without loss of generality. Call $\gamma \in \Gamma^+$ the curve where $f$ is discontinuous and let

$$K_1 = \{\mathbf{k} : \gamma \cap (B_{\mathbf{k}}^{h_1}/n) \neq \varnothing\}.$$

Here too, $|U_n^{h_1}(\mathbf{k})| \leq \beta h_1$ for $\mathbf{k} \notin K_1$ and $|U_n^{h_1}(\mathbf{k})| \leq 1$ for $\mathbf{k} \in K_1$. Also, $\#K_1 \leq Cn_1^2 h_1$, which comes from the fact that $\mathbf{k} \in K_1$ implies $\delta(\mathbf{k}) \leq \sqrt{2} h_1$ and the application of Lemma 16.1.

Using these facts and following the exact same arguments as for the one-dimensional case, we get

$$\mathcal{R}_n(M^{h_1,h_2}; f) \leq C \frac{1}{n_1^2} \sum_{\mathbf{k} \in \mathbb{I}_{n_1}^2} \mathbb{E}[(M^{h_2}[Y_n'](\mathbf{k}) - f(\mathbf{k}/n_1))^2] + Ch_1.$$

Using the upper bound (12.1) on the first term, we get

$$\mathcal{R}_n(M^{h_1,h_2}; f) \le C\left(h_2^2 + \frac{(\sigma_n')^2}{(n_1 h_2)^2} + h_2(\sigma_n')^2\right) + Ch_1.$$

Note that here $\sigma_n' \asymp (nh_1)^{-1}$. We then minimize over $h_1$ and $h_2$, with $h_1 = n^{-6/7}$ and $h_2 = n^{-4/7}$, and obtain the desired upper bound valid for any $f \in \mathcal{F}_2^+$.

14.2. *Lower bound.* The proof is completely parallel to the one-dimensional case, this time using (12.2).

## 15. Variability of quantiles.

LEMMA 15.1. *Assume $\Psi$ satisfies* [Shape] *and* [Decay]. *Then for all $\alpha_1 < \zeta/(\zeta - 1) < \alpha_2$, there are positive constants $C_1, C_2$ such that*

$$C_1(p(1-p))^{-\alpha_1} \le \frac{d}{dp}\Psi^{-1}(p) \le C_2(p(1-p))^{-\alpha_2} \qquad \forall p \in (0,1).$$

*Moreover, if $\psi(x)x^\zeta \asymp 1$, then*

$$\frac{d}{dp}\Psi^{-1}(p) \asymp (p(1-p))^{-\alpha},$$

*where $\alpha = \zeta/(\zeta - 1)$.*

PROOF. Let $1 < s < \zeta < t$ such that $\alpha_1 < s/(t-1) < t/(s-1) < \alpha_2$. We have

$$A_2(1+|x|)^{-t} \le \psi(x) \le A_1(1+|x|)^{-s} \qquad \forall x \in \mathbb{R}.$$

By integration, we also have

$$B_2(1+|x|)^{-t+1} \le 1 - \Psi(x) \le B_1(1+|x|)^{-s+1} \qquad \forall x \ge 0.$$

Therefore,

$$C_2^{-1}(1-\Psi(x))^{t/(s-1)} \le \psi(x) \le C_1^{-1}(1-\Psi(x))^{s/(t-1)} \qquad \forall x \ge 0.$$

By symmetry, we thus have

$$C_2^{-1}(\Psi(x)(1-\Psi(x)))^{t/(s-1)} \le \psi(x) \le C_1^{-1}(\Psi(x)(1-\Psi(x)))^{s/(t-1)} \qquad \forall x \in \mathbb{R}.$$

This is equivalent to

$$C_1(p(1-p))^{-s/(t-1)} \le \frac{d}{dp}\Psi^{-1}(p) \le C_2(p(1-p))^{-t/(s-1)} \qquad \forall p \in (0,1).$$

For the last statement, follow the same steps. □

LEMMA 15.2. *Let $\Psi$ satisfy* [Shape] *and* [Decay]. *Then for all $\alpha_1 < \zeta/(\zeta - 1) < \alpha_2$, there are positive constants $C_1, C_2$ such that*

$$C_1(p(1-p))^{-\alpha_1+1} \le \Psi^{-1}(p) \le C_2(p(1-p))^{-\alpha_2+1} \qquad \forall p \in (0,1).$$

*Moreover, if $\psi(x)x^\zeta \asymp 1$, then*

$$\Psi^{-1}(p) \asymp (p(1-p))^{-\alpha+1},$$

*where $\alpha = \zeta/(\zeta - 1)$.*

PROOF. Integrate the result in Lemma 15.1. □

15.1. *Proof of Lemma 6.2.* PROOF. Assume $\ell$ is odd for simplicity and let $\ell = 2m + 1$. We have

$$\Psi_{2m+1}(x) = (B_m \circ \Psi)(x/\sqrt{2m+1}),$$

where (see, e.g., [30]) $B_m$ is the $\beta$-distribution with parameters $(m, m)$:

$$B_m(y) = \frac{(2m+1)!}{(m!)^2} \int_0^y (u(1-u))^m \, du.$$

Given that $\Psi$ has a continuous density $\psi$ and $B_m$ is continuously differentiable, $\Psi_{2m+1}$ has a continuous density given by

$$\psi_{2m+1}(x) = \frac{1}{\sqrt{2m+1}} \psi(x/\sqrt{2m+1}) \cdot (B'_m \circ \Psi)(x/\sqrt{2m+1}).$$

Moreover, since $\psi$ is unimodal and symmetric about 0 and $B'_m$ is unimodal and symmetric about $1/2$, $\psi_{2m+1}$ is unimodal and symmetric about 0. Therefore, $\Psi_{2m+1}$ satisfies [Shape]. $\Psi_{2m+1}$ also satisfies [Decay] since $\psi_{2m+1}(x) \le C_m \psi(x/\sqrt{2m+1})$.

We now show that there is a constant $C$ such that, for $m$ large enough, $\psi_{2m+1}(x)(1 + |x|)^4 \le C$ for all $x$. It is enough to consider $x > 0$, which we do. Fix $s \in (1, \zeta)$. Using Stirling's formula and the fact that $\psi$ is bounded, we find $C$ such that

$$\psi_{2m+1}(x) \le C(4\Psi(x/\sqrt{2m+1})(1 - \Psi(x/\sqrt{2m+1})))^m.$$

In particular, $\psi_{2m+1}(x) \le C$ for all $x$. Since $\Psi(x)(1+x)^{s-1} \to 0$ as $x \to \infty$, there is $x_0 > 0$ such that $1 - \Psi(x) \le (1+x)^{-s+1}/4$ for $x \ge x_0$. Now, for $x \le x_0$, $(1+x)^4\psi_{2m+1}(x) \le C(1+x_0)^4$; for $x > x_0$,

$$(1+x)^4 \psi_{2m+1}(x) \le C \frac{(1+x)^4}{(1 + x/\sqrt{2m+1})^{(s-1)m}}.$$

By elementary calculus, as soon as $(s-1)m \ge 4\sqrt{2m+1}$, which happens when $m$ is large enough, the right-hand side is bounded by its value at $x_0$, which is also bounded by $C(1+x_0)^4$. □

**16. Some properties of planar curves.** This section borrows notation from Sections 10 and 11.

LEMMA 16.1.    *For* $\gamma \in \Gamma_2(\lambda)$ *and* $h \geq 1/n$, $\#\Delta \leq Cn^2h$, $C = 20(\lambda + 1)$.

PROOF.    Assume $\gamma$ is parametrized by arclength. Let $s_k = h/2 + kh$, for $k = 1, \ldots, [\text{length}(\gamma)/h]$. By the triangle inequality, for each $x \in [0,1]^2$ in the $h$-neighborhood of $\gamma$, there is $k = 1, \ldots, [\text{length}(\gamma)/h]$ such that $x$ and $\gamma(s_k)$ are within distance $2h$. Each ball centered at $\gamma(s_k)$ and of radius $2h$, $h \geq n^{-1}$, contains at most $20n^2h^2$ gridpoints. Therefore, the $h$-neighborhood of $\gamma$ contains at most $[\text{length}(\gamma)/h] \cdot Cn^2h^2 \leq C(\lambda+1)n^2h$ gridpoints.   □

LEMMA 16.2.    *There are constants* $h_0, C_1, C > 0$ *such that, if* $h < h_0$, *then for all* $\mathbf{i} \in \Delta$ *satisfying* $\delta(\mathbf{i}) > 2(C_1h^2 + n^{-1})$, $p(\mathbf{i}) \leq 1 - C\delta(\mathbf{i})/h$.

PROOF.    Let $h_0$ be defined as in Lemma 16.3 and assume $h < h_0$. Take $x$ so that $\gamma \cap B(x,h) \neq \varnothing$, where $B(x,h)$ is the disc of radius $h$ centered at $x$. By Lemma 16.3, there are arclengths $s_1 < s_2$ such that

$$\gamma \cap B(x,h) = \{\gamma(s) : s_1 < s < s_2\}.$$

A Taylor expansion of degree 2 gives

$$|\gamma(t) - \gamma(s) - (t-s)\gamma'(s)| \leq \kappa/2(t-s)^2 \qquad \forall s, t \in [0, \text{length}(\gamma)].$$

Together with the triangle inequality and the fact that $|\gamma'(s)| = 1$ for all $s$ and $|\gamma(s_2) - \gamma(s_1)| \leq 2h$, this implies

$$s_2 - s_1 \leq \kappa/2(s_2 - s_1)^2 + 2h.$$

Therefore, there is $C > 0$ such that $s_2 - s_1 \leq Ch$. Applying this Taylor expansion twice also implies

$$\left| \gamma(s) - \gamma(s_1) - (s - s_1)\frac{\gamma(s_2) - \gamma(s_1)}{s_2 - s_1} \right| \leq \kappa(s_2 - s_1)^2 \qquad \forall s \in [s_1, s_2],$$

which now becomes

$$\left| \gamma(s) - \gamma(s_1) - (s - s_1)\frac{\gamma(s_2) - \gamma(s_1)}{s_2 - s_1} \right| \leq C_1h^2 \qquad \forall s \in [s_1, s_2]$$

for some constant $C_1 > 0$. This means that, for all $s \in [s_1, s_2]$, $\gamma(s)$ is within distance $C_1h^2$ from the segment joining $\gamma(s_1)$ and $\gamma(s_2)$. Let $L$ be the line parallel to, and at distance $C_1h^2$ from $[\gamma(s_1), \gamma(s_2)]$, that is, closest to $x$. The line $L$ divides $B(x,h)$ into two parts $A$ and $B(x,h) \cap A^c$, where $x \in A$. Since we have $d(x, L) \geq d(x, \gamma) - d(L, \gamma) = d(x, \gamma) - C_1h^2$, if $d(x, \gamma) > C_1h^2$, $A \cap \gamma = \varnothing$ and $A$ contains the closed half disc with diameter parallel to $L$ that contains $x$.

Now, let $x$ be of the form $\mathbf{i}/n$ with $\mathbf{i} \in \mathbb{I}_n^2$ with $2C_1 h^2 + 2n^{-1} \leq \delta(\mathbf{i}) < h/2$. Without loss of generality, assume that $nh \geq 2$.

By symmetry, all open half discs of $B(\mathbf{i}/n, h)$ contain the same number of gridpoints, so that any closed half disc of $B(\mathbf{i}/n, h)$ contains more than half of the gridpoints within $B(\mathbf{i}/n, h)$. Let $A = A_0 \cup H$, where $H$ is the closed half disc with diameter parallel to $L$ that does not intersect $\gamma$. We have $\mathcal{G}(\mathbf{i}) = \mathbb{I}_n^2 \cap nA$, so that $\#\mathcal{G}(\mathbf{i}) \geq \#\mathcal{W}(\mathbf{i})/2 + \#(\mathbb{I}_n^2 \cap nA_0)$. $A_0$ contains a rectangular region $R$ with dimensions $d(\mathbf{i}/n, L)$ by $2\sqrt{h^2 - d(\mathbf{i}/n, L)^2}$, with $d(\mathbf{i}/n, L) = \delta(\mathbf{i}) - C_1 h^2 \geq 2/n$ and $2\sqrt{h^2 - d(\mathbf{i}/n, L)^2} \geq \sqrt{3}h \geq 2/n$. For such a rectangular region, with sidelengths of at least $2/n$,

$$R \subset \bigcup_{\mathbf{i}/n \in R} B(\mathbf{i}/n, 2/n),$$

so that

$$\#(\mathbb{I}_n^2 \cap nA_0) \geq \#(\mathbb{I}_n^2 \cap nR) \geq C|R|n^2 \geq Ch\delta(\mathbf{i})n^2.$$

It follows that $\rho(\mathbf{i}) \geq 1/2 + C\delta(\mathbf{i})/h$, which in turn implies $p(\mathbf{i}) \leq 1 - C\delta(\mathbf{i})/h$. We proved this for $\mathbf{i}$ such that $\delta(\mathbf{i}) < h/2$; however, this obviously extends to $\mathbf{i}$ such that $\delta(\mathbf{i}) < h$ with possibly a different constant $C$. $\quad\square$

LEMMA 16.3. *There is a constant $h_0 > 0$ such that the following holds for all $\gamma \in \Gamma^+$. If $h < h_0$ and $x \in [0,1]^2$ are such that $\gamma \cap B(x,h) \neq \varnothing$, then there are arclengths $s_1 < s_2$ such that $\gamma \cap B(x,h) = \{\gamma(s) : s_1 < s < s_2\}$.*

PROOF. We assume $\gamma$ is parametrized by arclength and consider arclengths modulo length($\gamma$).

Take $x \in [0,1]^2$ and $h > 0$ such that $\gamma \cap B(x,h) \neq \varnothing$. If $\gamma \subset \overline{B(x,h)}$, then $\gamma$ has maximum curvature bounded below by $h^{-1}$. We arrive at the same conclusion if $\gamma \cap \partial B(x,h)$ has infinite cardinality, for then $\gamma \cap \partial B(x,h)$ would have at least one accumulation point (since it is compact) at which the curvature would be exactly $h^{-1}$. Suppose $h < \kappa^{-1}$ so that $\gamma$ is not included in $\overline{B(x,h)}$ and $\gamma \cap \partial B(x,h)$ is nonempty and finite. Consider the set of arclengths $s$ with the property that there exists $\varepsilon_0 > 0$ such that, for all $0 < \varepsilon < \varepsilon_0$, $\gamma(s - \varepsilon) \in B(x,h)$ and $\gamma(s + \varepsilon) \notin B(x,h)$; this set is discrete, and therefore of the form $\{0 \leq s_1 < \cdots < s_m < \text{length}(\gamma)\}$. Note that because $\gamma$ is closed, $m$ is even.

Define $s_{m+1} = \text{length}(\gamma)$. We may assume that $s_1 = 0$ and $\gamma(s) \in \overline{B(x,h)}$ for all $s \in [s_1, s_2]$. Then, for all $k = 1, \ldots, m/2$, $\gamma(s) \notin \overline{B(x,h)}$ for all $s \in (s_{2k}, s_{2k+1})$.

Because $\gamma \in \Gamma^+$, we have, for all $k = 1, \ldots, m/2$,

$$\frac{s_{2k+1} - s_{2k}}{|\gamma(s_{2k+1}) - \gamma(s_{2k})|} \leq \theta \quad \text{and} \quad \frac{\text{length}(\gamma) - s_{2k+1} + s_{2k}}{|\gamma(s_{2k+1}) - \gamma(s_{2k})|} \leq \theta.$$

Suppose $m > 2$, so that $m \geq 4$. If $s_3 - s_2 \leq \text{length}(\gamma) - s_3 + s_2$, then $s_3 - s_2 \leq \theta|\gamma(s_3) - \gamma(s_2)| \leq 2\theta h$. Otherwise $s_3 - s_2 > \text{length}(\gamma) - s_3 + s_2$, and since $s_1 < s_2 < s_3 < s_m$, this implies that $\text{length}(\gamma) - s_m + s_1 \leq s_m - s_1$ and so $\text{length}(\gamma) - s_m + s_1 \leq \theta|\gamma(s_m) - \gamma(s_1)| \leq 2\theta h$. This is in turn equivalent to $s_{m+1} - s_m \leq \theta|\gamma(s_{m+1}) - \gamma(s_m)| \leq 2\theta h$. In both cases, there is $k = 1, \ldots, m/2$ such that $s_{2k+1} - s_{2k} \leq |\gamma(s_{2k+1}) - \gamma(s_{2k})| \leq 2\theta h$. Fix such a $k$. Let $a$ be the angle between $\gamma'(s_{2k})$ and $\gamma'(s_{2k+1})$ and let $b$ be the angle between $[x, \gamma(s_{2k})]$ and $[x, \gamma(s_{2k+1})]$. We have

$$\cos(a) = \langle \gamma'(s_{2k}), \gamma'(s_{2k+1}) \rangle = 1 - \frac{|\gamma'(s_{2k}) - \gamma'(s_{2k+1})|^2}{2}$$

with $|\gamma'(s_{2k}) - \gamma'(s_{2k+1})| \leq \kappa(s_{2k+1} - s_{2k}) \leq 2\kappa/\theta h$. Suppose $h < (\sqrt{2}\kappa\theta)^{-1}$, so that $a \leq C_1(s_{2k+1} - s_{2k})$, where $C_1 = C_1(\kappa, \theta)$. We also have

$$\sin(b/2) = \frac{|\gamma(s_{2k}) - \gamma(s_{2k+1})|}{2h} \geq \frac{s_{2k+1} - s_{2k}}{2\theta h},$$

so that $b \geq C_2(s_{2k+1} - s_{2k})/h$, where $C_2 = C_2(\kappa, \theta)$. Now, because $\gamma'(s_{2k})$ is either tangent or pointing outward and $\gamma'(s_{2k+1})$ is either tangent or pointing inward with respect to $B(x, h)$, we have $a \geq b$, if they are both tangent to $B(x, h)$, $a = b$. Therefore, $h \geq C_2/C_1$.

We thus let $h_0 = h_0(\kappa, \theta)$ be the minimum over all the constraints on $h$ and $C_2/C_1$. $\quad\square$

LEMMA 16.4. *There is a constant $C > 0$ such that, for all $h > 0$ and $\ell = 0, \ldots, nh$,*

$$\#\Xi_\ell \leq Cn.$$

PROOF. For $\mathbf{i} \in \Xi_\ell$, $B(\mathbf{i}/n, 1/(2n)) \subset T$, where $T = \{x : \ell - 1 \leq nd(x, \gamma) < \ell + 2\}$. Since those balls do not intersect, we have $\#\Xi_\ell C_1/n^2 \leq |T|$, $C_1 = \pi/4$.

As in the proof of Lemma 16.1, assume $\gamma$ is parametrized by arclength. Let $s_k = 1/(2n) + k/n$, for $k = 1, \ldots, [n \text{length}(\gamma)]$. Let $\vec{n}(s)$ be the normal vector to $\gamma$ at $\gamma(s)$ pointing out. Define $x_k^\pm = \gamma(s_k) \pm (\ell/n)\vec{n}(s_k)$. Take $x \in T$, say outside of $\gamma$; it is of the form $\gamma(s) + a\vec{n}(s)$, with $a \in [(\ell - 1)/n, (\ell + 2)/n]$. Let $k$ be such that $|s - s_k| \leq 1/n$. By the triangle inequality, we have

$$|x - x_k^+| \leq |\gamma(s) - \gamma(s_k)| + |a - \ell/n| + |\vec{n}(s) - \vec{n}(s_k)|$$

$$\leq |s - s_k|(1 + \kappa) + |a - \ell/n|$$

$$\leq C_2/n, \qquad C_2 = 3 + \kappa;$$

here we used $|\vec{n}(s) - \vec{n}(s_k)| = |\gamma'(s) - \gamma'(s_k)| \leq \kappa|s - s_k|$. Therefore, $T \subset \bigcup_k B(x_k^\pm, C_2/n)$, so that $|T| \leq n \cdot \text{length}(\gamma) \cdot C_2/n^2 \leq C_2 \cdot \text{length}(\gamma)/n$.

In the end, we have $\#\Xi_\ell \leq C_1|T|n^2 \leq C_1 \cdot C_2 \cdot \text{length}(\gamma) \cdot n$. $\quad\square$

# REFERENCES

[1] ANONYMOUS. (2007). Median filter. Wikipedia.

[2] BARNER, K. and ARCE, G. R. (2003). *Nonlinear Signal and Image Processing: Theory, Methods, and Applications.* CRC Press, Boca Raton, FL.

[3] BOTTEMA, M. J. (1991). Deterministic properties of analog median filters. *IEEE Trans. Inform. Theory* **37** 1629–1640. MR1134302

[4] BRANDT, J. (1998). Cycles of medians. *Util. Math.* **54** 111–126. MR1658177

[5] CANDÈS, E. J. and DONOHO, D. L. (2002). Recovering edges in ill-posed inverse problems: Optimality of curvelet frames. *Ann. Statist.* **30** 784–842. MR1922542

[6] CANDÈS, E. J. and DONOHO, D. L. (2002). Recovering edges in ill-posed inverse problems: Optimality of curvelet frames. *Ann. Statist.* **30** 784–842. MR1922542

[7] CAO, F. (1998). Partial differential equations and mathematical morphology. *J. Math. Pures Appl.* **77** 909–941. MR1656780

[8] CASELLES, V., SAPIRO, G. and CHUNG, D. H. (2000). Vector median filters, inf-sup operations, and coupled PDEs: Theoretical connections. *J. Math. Imaging Vision* **12** 109–119. MR1745601

[9] DONOHO, D. L. (1999). Wedgelets: Nearly minimax estimation of edges. *Ann. Statist.* **27** 859–897. MR1724034

[10] DONOHO, D. L. and JOHNSTONE, I. M. (1998). Minimax estimation via wavelet shrinkage. *Ann. Statist.* **26** 879–921. MR1635414

[11] DONOHO, D. L., JOHNSTONE, I. M., KERKYACHARIAN, G. and PICARD, D. (1995). Wavelet shrinkage: Asymptopia? *J. Roy. Statist. Soc. Ser. B* **57** 301–369. MR1323344

[12] DONOHO, D. L. and YU, T. P.-Y. (2000). Nonlinear pyramid transforms based on median-interpolation. *SIAM J. Math. Anal.* **31** 1030–1061. MR1759198

[13] FAN, J. and HALL, P. (1994). On curve estimation by minimizing mean absolute deviation and its implications. *Ann. Statist.* **22** 867–885. MR1292544

[14] GU, J., MENG, M., COOK, A. and FAULKNER, M. G. (2000). Analysis of eye tracking movements using fir median hybrid filters. In *ETRA '00: Proceedings of the 2000 symposium on Eye tracking research and applications* 65–69. ACM Press, New York.

[15] GUICHARD, F. and MOREL, J.-M. (1997). Partial differential equations and image iterative filtering. In *The State of the Art in Numerical Analysis (York, 1996). Inst. Math. Appl. Conf. Ser. New Ser.* **63** 525–562. Oxford Univ. Press, New York. MR1628359

[16] GUPTA, M. and CHEN, T. (2001). Vector color filter array demosaicing. In *Sensors and Camera Systems for Scientific, Industrial, and Digital Photography Applications.* II (M. B. J. C. N. Sampat, ed.). *Proceedings of the SPIE* **4306** 374–382. SPIE, Bellingham, WA.

[17] HAMZA, A. B., LUQUE-ESCAMILLA, P. L., MARTÍNEZ-AROZA, J. and ROMÁN-ROLDÁN, R. (1999). Removing noise and preserving details with relaxed median filters. *J. Math. Imaging Vision* **11** 161–177. MR1727352

[18] HUBER, P. J. (1964). Robust estimation of a location parameter. *Ann. Math. Statist.* **35** 73–101. MR0161415

[19] JUSTUSSON, B. (1981). Median filtering: Statistical properties. In *Two-Dimensional Digital Signal Processing. II* (T. S. Huang, ed.). *Topics in Applied Physics* **43** 161–196. Springer, Berlin. MR0688317

[20] KOCH, I. (1996). On the asymptotic performance of median smoothers in image analysis and nonparametric regression. *Ann. Statist.* **24** 1648–1666. MR1416654

[21] KOROSTELEV, A. P. and TSYBAKOV, A. B. (1993). *Minimax Theory of Image Reconstruction. Lecture Notes in Statistics* **82**. Springer, New York. MR1226450

[22] MALLOWS, C. L. (1979). Some theoretical results on Tukey's $3R$ smoother. In *Smoothing Techniques for Curve Estimation (Proc. Workshop, Heidelberg, 1979). Lecture Notes in Math.* **757** 77–90. Springer, Berlin. MR0564253

[23] MASON, D. M. (1984). Weak convergence of the weighted empirical quantile process in $L^2(0, 1)$. *Ann. Probab.* **12** 243–255. MR0723743

[24] MATHERON, G. (1975). *Random Sets and Integral Geometry.* Wiley, New York. MR0385969

[25] MOREL, J.-M. and SOLIMINI, S. (1995). *Variational Methods in Image Segmentation.* Birkhäuser Boston, Boston, MA. MR1321598

[26] MUMFORD, D. and SHAH, J. (1989). Optimal approximations by piecewise smooth functions and associated variational problems. *Comm. Pure Appl. Math.* **42** 577–685. MR0997568

[27] PERONA, P. and MALIK, J. (1990). Scale-space and edge detection using anisotropic diffusion. *IEEE Trans. Pattern Anal. Mach. Intell.* **12** 629–639.

[28] PITERBARG, L. I. (1984). Median filtering of random processes. *Problemy Peredachi Informatsii* **20** 65–73. MR0776767

[29] ROUSSEEUW, P. J. and BASSETT, G. W. JR. (1990). The remedian: A robust averaging method for large data sets. *J. Amer. Statist. Assoc.* **85** 97–104. MR1137355

[30] ROUSSEEUW, P. J. and BASSETT, G. W. JR. (1990). The remedian: A robust averaging method for large data sets. *J. Amer. Statist. Assoc.* **85** 97–104. MR1137355

[31] SAPIRO, G. (2001). *Geometric Partial Differential Equations and Image Analysis.* Cambridge Univ. Press, Cambridge. MR1813971

[32] SEMMES, S. W. (1988). Quasiconformal mappings and chord-arc curves. *Trans. Amer. Math. Soc.* **306** 233–263. MR0927689

[33] SERRA, J. (1982). *Image Analysis and Mathematical Morphology.* Academic Press, London. MR0753649

[34] SETHIAN, J. A. (1999). *Level Set Methods and Fast Marching Methods*, 2nd ed. *Cambridge Monographs on Applied and Computational Mathematics* **3**. Cambridge Univ. Press, Cambridge. MR1700751

[35] SHORACK, G. R. and WELLNER, J. A. (1986). *Empirical Processes with Aplications to Statistics.* Wiley, New York. MR0838963

[36] STRANNEBY, D. (2001). *Digital Signal Processing: DSP and Applications.* Oxford Univ. Press, London, UK.

[37] TUKEY, J. W. (1977). *Exploratory Data Analysis.* Addison-Wesley, Reading, MA.

[38] VELLEMAN, P. and HOAGLIN, D. (1981). *Applications, Basics, and Computing of Exploratory Data Analysis.* Duxbury, North Scituate, MA.

DEPARTMENT OF MATHEMATICS                    DEPARTMENT OF STATISTICS
UNIVERSITY OF CALIFORNIA, SAN DIEGO          STANFORD UNIVERSITY
9500 GILMAN DRIVE                            390 SERRA MALL
LA JOLLA, CALIFORNIA 92093-0112              STANFORD, CALIFORNIA 94305-4065
E-MAIL: eariasca@math.ucsd.edu               E-MAIL: donoho@stat.stanford.edu